\documentclass[12pt,a4paper]{article}
\usepackage{amsfonts,amsmath,amssymb,amsthm,accents,hyperref, xcolor}
\newcommand{\dis}{\displaystyle}
\newcommand{\txt}{\textstyle}

\newcommand{\noi}{\noindent}
\newcommand{\halmos}{\rule{1ex}{1.4ex}}
\newcommand{\QED}{\nopagebreak{\hspace*{\fill}$\halmos$\medskip}}

\newcommand{\med}{\medskip}
\newcommand{\quand}{\quad\mbox{and}\quad}

\newtheoremstyle{mythm}% name
  {}%      Space above
  {}%      Space below
  {\itshape}%         Body font
  {}%         Indent amount (empty = no indent, \parindent = para indent)
  {\bfseries}% Thm head font
  {}%        Punctuation after thm head
  {.5em}%     Space after thm head: " " = normal interword space;
  {#1 #2 \thmnote{(#3)}}%  Thm head spec (can be left empty, meaning `normal')

\theoremstyle{mythm}
\newtheorem{theorem}{Theorem}%[section]
\newtheorem{proposition}[theorem]{Proposition}
\newtheorem{lemma}[theorem]{Lemma}
\newtheorem{exercise}[theorem]{Exercise}
\newtheorem{corollary}[theorem]{Corollary}
\newtheorem{conjecture}[theorem]{Conjecture}

\newtheorem{counterex}[theorem]{Counterexample}

\newcommand{\bt}{\begin{theorem}}
\newcommand{\et}{\end{theorem}}
\newcommand{\bl}{\begin{lemma}}
\newcommand{\el}{\end{lemma}}
\newcommand{\bp}{\begin{proposition}}
\newcommand{\ep}{\end{proposition}}
\newcommand{\bcor}{\begin{corollary}}
\newcommand{\ecor}{\end{corollary}}
\newcommand{\br}{\begin{remark}\rm}
\newcommand{\er}{\end{remark}}
\newcommand{\bcon}{\begin{conjecture}}
\newcommand{\econ}{\end{conjecture}}
\newcommand{\bex}{\begin{exercise}}
\newcommand{\eex}{\end{exercise}}
\newcommand{\bcou}{\begin{counterex}}
\newcommand{\ecou}{\end{counterex}}

\newenvironment{Proof}[1][]{\noi\textbf{Proof #1}}{\QED}
\newcommand{\bpro}{\begin{Proof}}
\newcommand{\epro}{\end{Proof}}

\newcommand{\be}{\begin{equation}}
\newcommand{\ee}{\end{equation}}
\newcommand{\ba}{\begin{array}}
\newcommand{\bac}{\begin{array}{r@{\,}c@{\,}l}}
\newcommand{\ea}{\end{array}}
\newcommand{\bc}{\be\begin{array}{r@{\,}c@{\,}l}}
\newcommand{\ec}{\end{array}\ee}

\newcommand{\al}{\alpha}

\newcommand{\ga}{\gamma}
\newcommand{\Ga}{\Gamma}
\newcommand{\de}{\delta}
\newcommand{\De}{\Delta}
\newcommand{\eps}{\varepsilon}
\newcommand{\la}{\lambda}

\newcommand{\Om}{\Omega}

\newcommand{\Ci}{{\cal C}}

\newcommand{\Gi}{{\cal G}}
\newcommand{\Hi}{{\cal H}}
\newcommand{\Ii}{{\cal I}}

\newcommand{\Ki}{{\cal K}}
\newcommand{\Li}{{\cal L}}

\newcommand{\K}{{\mathbb K}}
\newcommand{\R}{{\mathbb R}}
\newcommand{\N}{{\mathbb N}}

\newcommand{\C}{{\mathbb C}}

\newcommand{\E}{{\mathbb E}}
\renewcommand{\P}{{\mathbb P}}

\newcommand{\li}{\langle}
\newcommand{\re}{\rangle}

\newcommand{\desd}{\ensuremath{\Leftrightarrow}}
\newcommand{\volgt}{\ensuremath{\Rightarrow}}

\newcommand{\sub}{\subset}
\newcommand{\beh}{\backslash}

\newcommand{\asto}[1]{\underset{{#1}\to\infty}{\longrightarrow}}

\newcommand{\ti}{\tilde}
\newcommand{\dgg}{\dagger}
\newcommand{\ov}{\overline}

\newcommand{\ffrac}[2]{{\textstyle\frac{{#1}}{{#2}}}}
\newcommand{\dif}[1]{\ffrac{\partial}{\partial{#1}}}
\newcommand{\diff}[1]{\ffrac{\partial^2}{{\partial{#1}}^2}}

\newcommand{\di}{\mathrm{d}}
\newcommand{\half}{{[0,\infty)}}
\newcommand{\expo}{\mbox{\large\it e}}
\newcommand{\ex}[1]{\expo^{\,\textstyle{#1}}}

\newcommand{\ha}{\ffrac{1}{2}}

\newcommand{\Xb}{{\mathbf X}}
\newcommand{\Yb}{{\mathbf Y}}

\newcommand{\tr}{{\rm tr}}
\newcommand{\ak}{\mathfrak{a}}
\newcommand{\bk}{\mathfrak{b}}
\newcommand{\ck}{\mathfrak{c}}
\newcommand{\gk}{\mathfrak{g}}
\newcommand{\hk}{\mathfrak{h}}
\newcommand{\ik}{\mathfrak{i}}
\newcommand{\jk}{\mathfrak{j}}

\newcommand{\ok}{\mathfrak{o}}

\newcommand{\uk}{\mathfrak{u}}
\newcommand{\glk}{\mathfrak{gl}}
\newcommand{\suk}{\mathfrak{su}}
\newcommand{\slk}{\mathfrak{sl}}
\newcommand{\sok}{\mathfrak{so}}

\newcommand{\ab}{\mathbf{a}}
\newcommand{\bb}{\mathbf{b}}
\newcommand{\cb}{\mathbf{c}}
\newcommand{\jb}{\mathbf{j}}
\newcommand{\kb}{\mathbf{k}}
\newcommand{\sv}{\mathbf{s}}
\newcommand{\tb}{\mathbf{t}}
\newcommand{\xb}{\mathbf{x}}
\newcommand{\yb}{\mathbf{y}}
\newcommand{\zb}{\mathbf{z}}

\begin{document}

%numbering formulas within sections
\makeatletter\@addtoreset{equation}{section}
\makeatother\def\theequation{\thesection.\arabic{equation}} 

\title{The Algebraic Approach to Duality:\\ An Introduction}
\author{Anja Sturm
\footnote{Institute for Mathematical Stochastics,
Georg-August-Universit\"at G\"ottingen,
Goldschmidtstr.~7,
37077 G\"ottingen,
Germany;
asturm@math.uni-goettingen.de}
\and
Jan M. Swart\footnote{The Czech Academy of Sciences,
Institute of Information Theory and Automation.
Pod vod\'arenskou v\v e\v z\' i 4,
18208 Prague 8,
Czech Republic;
swart@utia.cas.cz}
\and
Florian V\"ollering
\footnote{{Department of Mathematics, University of Bath, Claverton Down, Bath BA2 7AY, United Kingdom; f.m.vollering@bath.ac.uk}}}
\date{\today}
\maketitle

\begin{abstract}\noindent
This survey article gives an elementary introduction to the algebraic approach
to Markov process duality, as opposed to the pathwise approach. In the
algebraic approach, a Markov generator is written as the sum of products of
simpler operators, which each have a dual with respect to some duality
function. We discuss at length the recent suggestion by Giardin\`a, Redig, and
others, that it may be a good idea to choose these simpler operators in such a
way that they form an irreducible representation of some known Lie algebra.
In particular, we collect the necessary background on representations of Lie
algebras that is crucial for this approach. We also discuss older work by
Lloyd and Sudbury on duality functions of product form and the relation
between intertwining and duality.
\end{abstract}
\vspace{.5cm}

{\it Acknowledgement.} Work sponsored by grant 16-15238S of the Czech Science
Foundation (GA CR). %Kotecky-Swart 16-18

%82C22   	Interacting particle systems
%
%60K35   	Interacting random processes; statistical mechanics type
%               models; percolation theory
%17Bxx 	        Lie algebras and Lie superalgebras 
%17B10   	Representations, algebraic theory (weights)
%22E46   	Semisimple Lie groups and their representations
%
%22E60   	Lie algebras of Lie groups

\newpage

{\setlength{\parskip}{-2pt}\tableofcontents}

\newpage

\section{Introduction}\label{S:intro}

\subsection{Outline}\label{S:outline}

The aim of the present text is to give an introduction to the algebraic
approach to the theory of duality of Markov processes. In particular, we
present some of the pioneering work done by Lloyd and Sudbury
\cite{LS95,LS97,Sud00} and spend a lot of time explaining the more recent work
of Giardin\`a, Redig, and others \cite{GKRV09,CGGR15}. The algebraic approach
differs fundamentally from the pathwise approach propagated in e.g.,
\cite{JK14,SS16}. In principle, the algebraic approach is able to find a wider
class of dualities, but the price we pay for this is that it may suggest dual
operators that turn out not to be Markov generators.

In the remainder of this section, we quickly introduce the basic ideas behind
the algebraic approach. In Subsection~\ref{S:dual}, we explain how Markov
process duality can algebraically be viewed as an intertwining relation
between the generator of one Markov process and the adjoint of the generator
of another Markov process. As explained in Subsection~\ref{S:alg}, it is then
natural to view a Markov generator as being built up out of sums and
products of other, simpler operators. If all these building blocks have duals
with respect to a duality function, then so has the original Markov generator.

A central idea of of Giardin\`a, Redig, et al.\ \cite{GKRV09,CGGR15} is to
choose these building blocks so that they form a representation of a Lie
algebra. To understand why that may be a good idea, one needs quite a bit of
background on Lie algebras. Since probabilists may not be familiar with this,
after a small detour to pathwise duality in Subsection~\ref{S:path}, we devote
all of Section~\ref{S:Liealg} to providing this background.

The study of Lie algebras and their representations is a huge subject with a
venerable history. Although there exist good introductory texts, we will need
some theory that is considered too advanced for the usual textbooks.  In
particular, this refers to the representation theory of non-compact Lie groups
like SU(1,1) or the Heisenberg group. In order to squeeze the essential facts
that we need for our purposes into little over 10 pages, we had to cut some
corners and in some cases resign on full mathematical rigour. We also leave
out a lot of background material (e.g., Lie groups, as opposed to Lie
algebras, stay almost completely out of the picture). To partly compensate for
this, we have added Appendix~\ref{A:Lie} which gives a somewhat more complete,
but still sketchy picture.

After our little excursion into Lie algebras, in Section~\ref{S:repdu}, we
come to the core of our text. In Subsections~\ref{S:WFself}, \ref{S:SEP}, and
\ref{S:SIP} we demonstrate the approach via Lie algebras on three examples,
which are based on representation theory for the Heisenberg algebra, SU(2),
and SU(1,1), respectively. In Subsection~\ref{S:WFself}, we formulate a
general principle and apply it to discover a self-duality of the Wright-Fisher
diffusion from (\ref{WF}). After Subsection~\ref{S:contwine}, which is needed
to deal with infinite state spaces, in Subsection~\ref{S:SEP}, we use the
well-known representation theory of SU(2) to derive a duality for the
symmetric exclusion process. This duality is not very interesting on its own,
but serves as a preparation for the symmetric inclusion process in
Section~\ref{S:SIP} which turns out to be very similar to the former, except
that SU(2) is replaced by SU(1,1).

In Sections~\ref{S:sudbur}--\ref{S:bias} we present results of Lloyd and
Sudbury \cite{LS95,LS97,Sud00} that do not require knowledge of Lie algebras,
but do use some facts about tensor products from Section~\ref{S:sumprod}. In
particular, in Section~\ref{S:sudbur} we discuss duality functions of product
form, including q-duality, while in Section~\ref{S:twine} we discuss
intertwining of Markov processes, and in particular thinning relations which
are closely connected to q-duality.

In Sections~\ref{S:symm} and \ref{S:SEP2}, finally, we discuss another
observation from \cite{GKRV09}, who show that nontrivial dualities can
sometimes be found by starting from a ``trivial'' duality which is based on
reversibility, and then using a symmetry of the model to transform such a
duality into a nontrivial one. Although Lie algebras are not strictly needed
in this approach, writing generators in terms of the basis elements of a
representation of a Lie algebra can help finding suitable symmetries.

\subsection{Markov duality and intertwining}\label{S:dual}

In Section~\ref{S:intro}, for technical simplicity, we mostly restrict
ourselves to Markov processes with finite state spaces. As we will see in
Section~\ref{S:repdu}, many of the basic ideas discussed here can with some
care be made to work also in infinite dimensional settings. How to do this is
in part discussed in Section~\ref{S:contwine}, but for brevity, we will not
always go into the technical details and sometimes use the calculations of the
present section merely as an inspiration.

The generator of a continuous-time Markov process with finite state space
$\Om$ is a matrix $L$ such that
\be\label{gendef}
L(x,y)\geq 0\quad(x\neq y)\quand\sum_yL(x,y)=0.
\ee
Equivalently, we can identify $L$ with the linear operator
  $L:\R^\Om\to\R^\Om$ defined by
\be\label{Lin}
Lf(x):=\sum_{y\in\Om}L(x,y)f(y)\qquad(x\in\Om).
\ee
A linear operator $L:\R^\Om\to\R^\Om$ is a Markov generator (i.e., satisfies
(\ref{gendef})) if and only if the semigroup\footnote{The semigroup property
  says that $P_0=I$ and $P_sP_t=P_{s+t}$.} of operators $(P_t)_{t\geq 0}$
defined by
\[
P_t:=e^{tL}=\sum_{n=0}^\infty\frac{1}{n!}t^nL^n
\]
is a Markov semigroup, i.e., $P_t$ is a probability kernel for each $t\geq 0$.
If $L$ is a Markov generator, then $(P_t)_{t\geq 0}$ are the transition
kernels of some $\Om$-valued Markov process $(X_t)_{t\geq 0}$.

Let $\Om$ and $\hat\Om$ be finite sets. We can view a function
$D:\Om\times\hat\Om\to\R$ as a matrix
\[
(D(x,y))_{x\in\Om,\ y\in\hat\Om}
\]
that as in (\ref{Lin}) corresponds to a linear operator
$D:\R^{\hat\Om}\to\R^\Om$.

Let $L$ and $\hat L$ be generators of Markov processes $(X_t)_{t\geq 0}$ and
$(Y_t)_{t\geq 0}$ with state spaces $\Om$ and $\hat\Om$ and semigroups
$(P_t)_{t\geq 0}$ and $(\hat P_t)_{t\geq 0}$, and let
$D:\Om\times\hat\Om\to\R$ be a function. We make the following simple
observation. Below, we let $A^\dgg(x,y):=A(y,x)$ (or
$A^\dgg(x,y):=\ov{A(y,x)}$ for matrices over the complex numbers) denote the
\emph{adjoint} of a matrix $A$.\footnote{In other words, $\li A^\dgg
  f|g\re:=\li f|Ag\re$ where $\li f|g\re:=\sum_{x\in\Om}\ov{f(x)}g(x)$ denotes
  the usual inner product on $\C^\Om$.  For adjoints with respect to a general
  inner product on finite or infinite dimensional spaces we write $A^\ast$.}

\bl[Duality]
The\label{L:du} following conditions are equivalent.
\begin{enumerate}
\item $\dis LD=D\hat L^\dgg$,
\item $\dis P_tD=D\hat P^\dgg_t$ for all $t\geq 0$,
\item $\dis\E^x[D(X_t,y)]=\E^y[D(x,Y_t)]$ for all $x\in\Om$, $y\in\hat\Om$,
  and $t\geq 0$.
\end{enumerate}
\el
\bpro
If (i) holds for $L$, then it also holds for any linear combination of powers
of $L$. In particular, filling in the definition of $P_t$, we see that (i)
implies (ii). Conversely, differentiating with respect to $t$, we see that
(ii) implies (i). Condition~(iii) is just a rewrite of (ii).
\epro

If the conditions of Lemma~\ref{L:du} are satisfied, then we say that
$(X_t)_{t\geq 0}$ and $(Y_t)_{t\geq 0}$ are \emph{dual} with \emph{duality
  function} $D$. If $L=\hat L$, then we speak of \emph{self-duality}.
Condition~(i) can also be written as
\be\label{LDL}
LD(\,\cdot\,,y)(x)=\hat LD(x,\,\cdot\,)(y)\qquad(x\in\Om,\ y\in\hat\Om).
\ee
Under suitable assumptions, the equivalence of (iii) and (\ref{LDL}) can
often also be established for Markov processes with infinite state space.

An algebraic relation of the form $AB=BC$ is called an \emph{intertwining
  relation} between operators $A$ and $C$. The operator $B$ is called the
\emph{intertwiner}. Thus, Lemma~\ref{L:du} says that two Markov processes are
dual if and only if there exists an intertwiner between the generator of one
Markov process, and the adjoint of the generator of another Markov process.
Note that if $L$ is dual to $\hat L$ with duality function $D$, then $\hat L$
is dual to $L$ with duality function $D^\dgg$. Thus, duality is a symmetric
concept.

Closely related to Markov process duality is the concept of \emph{intertwining
  of Markov processes}, which has a more narrow meaning than the algebraic
concept of intertwining. Let, again, $L$ and $\hat L$ be generators of Markov
processes $(X_t)_{t\geq 0}$ and $(Y_t)_{t\geq 0}$ with state spaces $\Om$ and
$\hat\Om$ and semigroups $(P_t)_{t\geq 0}$ and $(\hat P_t)_{t\geq 0}$. Let
$K:\Om\times\hat\Om\to\R$ be a function. In what follows, we assume that $K$
is a probability kernel, i.e., $K(x,y)\geq 0$ $\forall x,y$ and
$\sum_yK(x,y)=1$ for each $x$.

\bl[Intertwining of Markov processes]
The\label{L:tw} the following conditions are equivalent.
\begin{enumerate}
\item $\dis LK=K\hat L$.
\item $\dis P_tK=K\hat P_t\qquad(t\geq 0)$.
\item $\dis\mu_0K=\nu_0$ implies $\dis\mu_0P_tK=\nu_0\hat P_t\qquad(t\geq 0)$.
\end{enumerate}
\el
\bpro
The equivalence of (i) and (ii) follows by the same argument as in
Lemma~\ref{L:du}. Condition (ii) implies $\mu_0P_tK=(\mu_0K)\hat P_t$ $(t\geq
0)$. Setting $\mu_0=\de_x$ we see that (iii) implies (ii).
\epro

In condition~(iii), note that $\mu_0P_t$ and $\nu_0\hat P_t$ describe the laws
at time $t$ of the Markov processes $(X_t)_{t\geq 0}$ and $(Y_t)_{t\geq 0}$
started in initial laws $\mu_0$ and $\nu_0$, respectively. If the conditions
of Lemma~\ref{L:tw} are satisfied, then we say that the Markov processes
$(X_t)_{t\geq 0}$ and $(Y_t)_{t\geq 0}$ are \emph{intertwined}.

If $K$ is invertible as a matrix, then $LK=K\hat L$ implies $\hat
LK^{-1}=K^{-1}L$; however, $K^{-1}$ will in general not be a probability
kernel. In view of this, in an intertwining relation between Markov processes,
the two processes do not play symmetric roles. To stress the different roles
of $X$ and $Y$, following \cite{Swa13}, it is convenient to say that $Y$ is an
intertwined Markov process \emph{on top} of $X$.

If the conditions of Lemma~\ref{L:tw} are satisfied, then the Markov processes
$X$ and $Y$ can actually be coupled such that $(X_t,Y_t)_{t\geq 0}$ is a
Markov process and
\[
\P[Y_t\in\,\cdot\,|(X_s)_{0\leq s\leq t}]=K(X_s,\,\cdot\,)
\quad{\rm a.s.}\quad(t\geq 0),
\]
see \cite{Fil92,Swa13}. Note that this strengthens condition~(iii) of
Lemma~\ref{L:tw}.

\subsection{The algebraic approach}\label{S:alg}

We make the following simple observation. Below, $\R^\Om$ denotes the space of
all functions $f:\Om\to\R$.

\bl[Duality of building blocks]
Let\label{L:build} $\Om,\hat\Om$ be finite spaces and let $A_i:\R^\Om\to\R^\Om$,
$B_i:\R^{\hat\Om}\to\R^{\hat\Om}$ $(i=1,2)$, and $D:\R^{\hat\Om}\to\R^\Om$ be
linear operators such that
\be
A_iD=DB_i^\dgg\qquad(i=1,2).
\ee
Then
\be
(r_1A_1+r_2A_2)D=D(r_1B_1+r_2B_2)^\dgg
\quand
(A_1A_2)D=D(B_2B_1)^\dgg.
\ee
\el

Lemma~\ref{L:build} implies that if we can write a Markov generator $L$ as a
linear combination of products of ``simpler'' operators $A_i$, for example,
(denoting the identity operator by $I$),
\be\label{building}
L=r_\emptyset I+r_1A_1+r_{23}A_2A_3+r_{113}A_1^2A_3,
\ee
and these ``building blocks'' satisfy $A_iD=DB_i^\dgg$ for some duality
function $D$, the $L$ will be dual to the operator
\be
\hat L=r_\emptyset I+r_1B_1+r_{23}B_3B_2+r_{113}B_3B_1^2.
\ee
Note that in each term, we have not only replaced $A_i$ by $B_i$ but also
reversed the order of the factors. If we are lucky, $\hat L$ is a Markov
generator and we have discovered a Markov duality.

We demonstrate this approach on the Wright-Fisher diffusion with selection
parameter $s\in\R$, which is the diffusion in $[0,1]$ with generator
\be\label{WF}
Lf(x)=x(1-x)\diff{x}+sx(1-x)\dif{x}.
\ee
We are immediately cheating here, since $L$ is not a linear operator acting on
a finite dimensional space. Ignoring the difficulties associated with infinite
dimension, we can write $L$ in terms of simpler ``building blocks'' as
follows. We set
\be\label{Apm}
A^-f(x):=(1-x)f(x)\quand A^+f(x):=\dif{x}f(x),
\ee
and we write $L$ in terms of these building blocks as
\be\label{LAexpr}
L=A^-(I-A^-)A^+(sI+A^+).
\ee
As our dual space, we choose $\N=\{0,1,\ldots\}$ and as our duality function
we choose the function $D:[0,1]\times\N\to\R$ given by
\be\label{Dmom}
D(x,n):=(1-x)^n\qquad(x\in[0,1],\ n\in\N).
\ee
Let $B^\pm$ be operators acting on functions $f:\N\to\R$ as
\be
B^-f(n):=f(n+1)\quand B^+f(n):=-nf(n-1).
\ee
Then $B^\pm$ are dual to $A^\pm$ in the sense of (\ref{LDL}), i.e.,
\be
A^\pm D(\,\cdot\,,n)(x)=B^\pm D(x,\,\cdot\,)(n)\qquad(x\in[0,1],\ n\in\N).
\ee
Therefore, in view of Lemma~\ref{L:build}, the following operator should be
dual to $L$:
\be
\hat L=(sI+B^+)B^+(I-B^-)B^-.
\ee
(Note that we have replaced $A^\pm$ by $B^\pm$ and reversed the order of the
factors.) A little calculation reveals that
\be\label{WFmomdu}
\hat Lf(n)=n(n-1)\big\{f(n-1)-f(n)\}+sn\big\{f(n+1)-f(n)\big\}.
\ee
This is not, in general, a Markov generator. For $s\geq 0$, however, it is the
generator of a Markov process in $\N$ that jumps from $n$ to $n-1$ with rate
$n(n-1)$ and from $n$ to $n+1$ with rate $sn$.

Recall that the \emph{commutator} of two operators $A,B$ is defined as
$[A,B]:=AB-BA$. For our operators $A^\pm$, it is easy to check that
\be\label{ApmHeis}
[A^-,A^+]=I.
\ee
This is similar to the commutation relation between the position and momentum
operators in quantum physics. Indeed, the operators $A^\pm$ can be used to
define a representation of the Heisenberg algebra, which is a particular Lie
algebra. The connection to Lie algebras can help us to choose good building
blocks and can sometimes also suggest duality functions. To explain this, we
need some theory about representations of Lie algebras, which will be
presented in the next section.

\subsection{The pathwise approach}\label{S:path}

In the remainder of this section, we point out some differences and
similarities between the algebraic and pathwise approaches to Markov process
duality. A \emph{random mapping representation} of a probability kernel $K$ is
a random map $M$ such that
\be
K(x,\di y)=\P[M(x)\in\di y].
\ee
A \emph{stochastic flow} is a collection $(\Xb_{s,u})_{s\leq u}$ of random
maps $\Xb_{s,u}:\Om\to\Om$ such that $\Xb_{s,s}=I$ and
$\Xb_{t,u}\circ\Xb_{s,t}=\Xb_{s,u}$. We say that $(\Xb_{s,u})_{s\leq u}$ has
\emph{independent increments} if
\be
\Xb_{t_1,t_2},\ldots,\Xb_{t_{n-1},t_n}
\ee
are independent for any $t_1<\cdots<t_n$. If $(\Xb_{s,u})_{s\leq u}$ is a
stochastic flow with independent increments such that the law of $\Xb_{s,u}$
depends only on the difference $u-s$, and $X_0$ is an independent
$\Om$-valued random variable, then setting
\be
X_t:=\Xb_{0,t}(X_0)\qquad(t\geq 0)
\ee
defines a Markov process with transition kernels
\be
P_{u-s}(x,\di y):=\P[\Xb_{s,u}(x)\in\di y]\qquad(s\leq u).
\ee
Note that this formula says that $\Xb_{s,u}$ is a random mapping
representation of $P_{u-s}$.

Markov processes can often be constructed from stochastic flows. For example,
if a stochastic differential equation has unique strong solutions, then these
solutions (for different initial states) define a stochastic flow with
independent increments that can be used to construct a diffusion process.
If $L$ is the generator of a Markov process with finite state space $\Om$,
then $L$ can always be written in the form
\be\label{Lrandmap}
Lf(x)=\sum_{m\in\Gi}r_m\big\{f\big(m(x)\big)-f\big(x\big)\big\},
\ee
where $\Gi$ is a finite collection of maps $m:\Om\to\Om$. We say that two maps
$m,\hat m$ are \emph{dual} with respect to a duality function $D$ if
\be \label{eq:pathD}
D\big(m(x),y\big)=D\big(x,\hat m(y)\big)\qquad(x\in\Om,\ y\in\hat\Om).
\ee
Two stochastic flows $(\Xb_{s,u})_{s\leq u}$ and $(\Yb_{s,u})_{s\leq u}$ are
dual\footnote{The definition of duality for stochastic flows that we give here
  is a weak one. It is often natural to give a somewhat stronger definition,
  see \cite{SS16}.} if for each $s\leq u$, a.s., $\Yb_{-u,-s}$ is dual to
$\Xb_{s,u}$. If two stochastic flows are dual, then we say that their
associated Markov processes are \emph{pathwise dual}. It is easy to see that
this implies Markov process duality.

We recall that in the algebraic approach, there may be many ways in which a
given Markov generator can be written in terms of more elementary ``building
blocks'' as in (\ref{building}). Similarly, in the pathwise approach, there
are usually many different ways in which a Markov generator can be written in
terms of maps as in (\ref{Lrandmap}).  In the algebraic approach we have seen
that if all building blocks have duals with respect to a given duality
function, then a Markov generator built up from these building blocks also has
a dual $\hat L$. Similarly, in the pathwise approach, if all maps $m$
occurring in (\ref{Lrandmap}) have duals $\hat m$ with respect to some duality
function $D$, then the process with generator $L$ is pathwise dual to the
process with generator
\be
\hat Lf(x):=\sum_{m\in\Gi}r_m\big\{f\big(\hat m(x)\big)-f\big(x\big)\big\}.
\ee
An advantage of the pathwise approach is that an operator $\hat L$ of this
form is guaranteed to me a Markov generator. On the other hand, not all
dualities can be constructed as pathwise dualities, so the algebraic approach
is more general. Nevertheless, many known dualities, including the duality for
the Wright-Fisher diffusion discussed in the previous subsection, can be
obtained in a pathwise way or as limits of such pathwise dualities, see
\cite{Swa06,AH07}.

There are more analogies between the algebraic and pathwise approaches. In
Subsection~\ref{S:symm}, we will see that in the algebraic approach,
nontrivial dualities can sometimes be found by starting with a ``trivial''
duality obtained from reversibility and then applying a symmetry
transformation. In \cite{SS16}, it is shown that nontrivial pathwise dualities
can be found by starting with a ``trivial'' duality to the inverse image map
and then looking for invariant subspaces of the dual process.

\section{Representations of Lie algebras}\label{S:Liealg}

\subsection{Lie algebras}\label{S:Lieint}

A complex\footnote{In this section, we mostly focus on complex Lie
  algebras. Some results stated in the present section (in particular,
  part~(b) of Schur's lemma) are true for complex Lie algebras only. See
  Appendix~\ref{A:Lie} for a more detailed discussion.} (resp.\ real)
\emph{Lie algebra} is a finite-dimensional linear space $\gk$ over $\C$
(resp.\ $\R$) together with a map $[\,\cdot\,,\,\cdot\,]:\gk\times\gk\to\gk$
called \emph{Lie bracket} such that
\begin{enumerate}
\item $(\xb,\yb)\mapsto[\xb,\yb]$ is bilinear,
\item $[\xb,\yb]=-[\yb,\xb]$ (skew symmetry),
\item $[\xb,[\yb,\zb]]+[\yb,[\zb,\xb]]+[\zb,[\xb,\yb]]=0$ (Jacobi identity).
\end{enumerate}
An \emph{adjoint operation} on a Lie algebra $\gk$ is a map
$\xb\mapsto\xb^\ast$ such that
\begin{enumerate}
\item $\xb\mapsto\xb^\ast$ is conjugate linear,
\item $(\xb^\ast)^\ast=\xb$,
\item $[\xb^\ast,\yb^\ast]=[\yb,\xb]^\ast$.
\end{enumerate}
If $\gk$ is a complex Lie algebra, then the space of its skew symmetric
elements $\hk:=\{\xb\in\gk:\xb^\ast=-\xb\}$ forms a real Lie algebra.
Conversely, starting from a real Lie algebra $\hk$, we can always find a
complex Lie algebra $\gk$ equipped with a adjoint operation such that $\hk$
is the space of skew symmetric elements of $\gk$. Then $\gk$ is called the
\emph{complexification} of $\hk$.

If $\{\xb_1,\ldots,\xb_n\}$ is a basis for $\gk$, then the Lie bracket on $\gk$ 
is uniquely characterized by the \emph{commutation relations}
\[
[\xb_i,\xb_j]=\sum_{k=1}^nc_{ijk}\xb_k\qquad(i<j).
\]
The constants $c_{ijk}$ are called the \emph{structure constants}. If $\gk$ is
equipped with an adjoint operation, then the latter is uniquely characterized
by the \emph{adjoint relations}
\[
\xb_i^\ast=\sum_jd_{ij}\xb_j.
\]

\noi
\textbf{Example} Let $V$ be a finite dimensional complex linear space, let
$\Li(V)$ denote the space of all linear operators $A:V\to V$, and let $\tr(A)$
denote the trace of an operator $A$. Then
\[
\gk:=\{A\in\Li(V):\tr(A)=0\}
\quad\mbox{with}\quad
[A,B]:=AB-BA
\]
is a Lie algebra. Note that $\tr([A,B])=\tr(AB)-\tr(BA)=0$ by
the basic property of the trace, which shows that $[A,B]\in\gk$ for all
$A,B\in\gk$. Note also that $\gk$ is in general not an algebra, i.e.,
$A,B\in\gk$ does not imply $AB\in\gk$. If $V$ is equipped with an inner
product $\li\,\cdot\,|\,\cdot\,\re$ (which we always take conjugate linear in its
first argument and linear in its second argument) and $A^\ast$ denotes the
adjoint of $A$ with respect to this inner product, i.e.,
\[
\li A^\ast v|w\re:=\li v|Aw\re,
\]
then one can check that $A\mapsto A^\ast$ is an adjoint operation on $\gk$.\med

By definition, a \emph{Lie algebra homomorphism} is a map
$\phi:\gk\to\hk$ from one Lie algebra into another that preserves the
structure of the Lie algebra, i.e., $\phi$ is linear and
\[
\phi([A,B])=[\phi(A),\phi(B)].
\]
If $\phi$ is invertible, then its inverse is also a Lie algebra
homomorphism. In this case we call $\phi$ a \emph{Lie algebra isomorphism}.
We say that a Lie algebra homomorphism $\phi$ is \emph{unitary} if it moreover
preserves the structure of the adjoint operation, i.e.,
\[
\phi(A^\ast)=\phi(A)^\ast.
\]

If $\gk$ is a Lie algebra, then we can define a \emph{conjugate} of $\gk$,
which is a Lie algebra $\ov\gk$ together with a conjugate linear bijection
$\gk\ni\xb\mapsto\ov\xb\in\ov\gk$ such that
\[
[\ov\xb,\ov\yb]=[\yb,\xb].
\]
It is easy to see that such a conjugate Lie algebra is unique up to natural
isomorphisms, and that the $\ov{\ov\gk}$ is naturally isomorphic to $\gk$.
If $\gk$ is equipped with an adjoint operation, then we can define an adjoint
operation on $\ov\gk$ by ${\ov\xb}^\ast:=\ov{(\xb^\ast)}$.\med

\noi
\textbf{Example} Let $V$ be a complex linear space on which an inner product
is defined and let $\gk\sub\Li(V)$ be a linear subspace such that $A,B\in\gk$
implies $[A,B]\in\gk$. Then $\gk$ is a sub-Lie-algebra of $\Li(V)$.
Now $\ov\gk:=\{A^\ast:A\in\gk\}$, together with the map $\ov A:=A^\ast$ is a
realization of the conjugate Lie algebra of $\gk$.

\subsection{Representations}\label{S:rep}

If $V$ is a finite dimensional linear space, then the space $\Li(V)$ of linear
operators $A:V\to V$, equipped with the \emph{commutator}
\[
[A,B]:=AB-BA
\]
is a Lie algebra. By definition, a \emph{representation} of a complex Lie
algebra $\gk$ is a pair $(V,\pi)$ where $V$ is a complex linear space of
dimension ${\rm dim}(V)\geq 1$ and $\pi:\gk\to\Li(V)$ is a Lie algebra
homomorphism. A representation is \emph{unitary} if this homomorphism is
unitary and \emph{faithful} if $\pi$ is an isomorphism to its image
$\pi(\gk):=\{\pi(\xb):\xb\in\gk\}$.

There is another way of looking at representations that is often useful. If
$(V,\pi)$ is a representation, then we can define a map
\[
\gk\times V\ni(\xb,v)\mapsto\xb v\in V
\]
by $\xb v:=\pi(\xb)v$. Such a map satisfies
\begin{enumerate}
\item $(\xb,v)\mapsto Av$ is bilinear (i.e., linear in both arguments),
\item $[\xb,\yb]v=\xb(\yb v)-\yb(\xb v)$.
\end{enumerate}
Any map with these properties is called a \emph{left action} of $\gk$ on $V$.
It is easy to see that if $V$ is a complex linear space that is equipped with
a left action of $\gk$, then setting $\pi(\xb)v:=\xb v$ defines a Lie algebra
homomorphism from $\gk$ to $\Li(V)$. Thus, we can view representations as
linear spaces on which a left action of $\gk$ is defined.\med

\noi
\textbf{Example} For any Lie algebra, we may set $V:=\gk$. Then, using the
Jacobi identity, one can verify that the map $(\xb,\yb)\mapsto[\xb,\yb]$ is a
left action of $\gk$ on $V$. (See Lemma~\ref{L:Lierep} in the appendix.)
In this way, every Lie algebra can be represented
on itself. This representation is not always faithful, but for many Lie
algebras of interest, it is.\med

Yet another way to look at representations is in terms of commutation
relations. Let $\gk$ be a Lie algebra with basis elements
$\xb_1,\ldots,\xb_n$, which satisfy the commutation relations
\[
[\xb_i,\xb_j]=\sum_{k=1}^nc_{ijk}\xb_k\qquad(i<j).
\]
Let $V$ be a complex linear space with ${\rm dim}(V)\geq 1$ and let
$X_1,\ldots,X_n\in\Li(V)$ satisfy
\[
[X_i,X_j]=\sum_{k=1}^nc_{ijk}X_k\qquad(i<j).
\]
Then there exists a unique Lie algebra homomorphism $\pi:\gk\to\Li(V)$ such
that $\pi(\xb_i)=X_i$ $(i=1,\ldots,n)$. Thus, any collection of linear
operators that satisfies the commutation relations of $\gk$ defines a
representation of $\gk$. Such a representation is faithful if and only if
$X_1,\ldots,X_n$ are linearly independent. If $\gk$ is equipped with an
adjoint operation and $V$ is equipped with an inner product, then the
representation $(V,\pi)$ is unitary if and only if
$X_1,\ldots,X_n$ satisfy the adjoint relations of $\gk$, i.e.,
\[
\xb_i^\ast=\sum_jd_{ij}\xb_j
\quand
X_i^\ast=\sum_jd_{ij}X_j.
\]

Let $V$ be a representation of a Lie algebra $\gk$. By definition, an
\emph{invariant subspace} of $V$ is a linear subspace $W\sub V$ such that $\xb
w\in W$ for all $w\in W$ and $\xb\in\gk$. A representation is
\emph{irreducible} if its only invariant subspaces are $W=\{0\}$ and $W=V$.

Let $V,W$ be two representations of the same Lie algebra $\gk$. By definition,
an \emph{intertwiner} of representations is a linear map $\phi:V\to W$ that
preserves the structure of a representation, i.e.,
\[
\phi(\xb v)=\xb\phi(v).
\]
If $\phi$ is a bijection then its inverse is also an intertwiner. In this
case we call $\phi$ an \emph{isomorphism} and say that the representations are
\emph{equivalent} (or \emph{isomorphic}).

The following result can be found in, e.g., \cite[Thm~4.29]{Hal03}.
Below and in what follows, we let $I\in\Li(V)$ denote the \emph{identity
  operator} $Iv:=v$.

\bp[Schur's lemma]\label{P:Schur}\quad
\begin{itemize}
\item[{\rm\textbf{(a)}}] Let $V$ and $W$ be irreducible representations of the
  same Lie algebra and let $\phi:V\to W$ be an intertwiner. Then either
  $\phi=0$ or $\phi$ is an isomorphism.
\item[{\rm\textbf{(b)}}] Let $V$ be an irreducible representation of a
  Lie algebra and let $\phi:V\to V$ be an intertwiner. Then $\phi=\la I$ for
  some $\la\in\C$.
\end{itemize}
\ep

For us, the following simple consequence of Schur's lemma will be important.

\bcor[Unique intertwiner]
Let\label{C:twin} $(V,\pi_V)$ and $(W,\pi_W)$ be equivalent irreducible
representations of some Lie algebra. Then there exists an intertwiner
$\phi:V\to W$ that is unique up to a multiplicative constant, such that
\[
\phi\pi_V(\xb)=\pi_W(\xb)\phi.
\]
\ecor
\bpro
By assumption, $V$ and $W$ are equivalent, so there exists an isomorphism
$\phi:V\to W$. Assume that $\psi:V\to W$ is another intertwiner. Then
$\phi^{-1}\circ\psi$ is an intertwiner from $V$ into itself, so by part~(b) of
Schur's lemma, $\phi^{-1}\circ\psi=\la I$ and hence $\psi=\la\phi$.
\epro

If $V$ is a complex linear space, then we can define a \emph{conjugate} of
$V$, which is a complex linear space $\ov V$ together with a conjugate linear
bijection $\phi\mapsto\ov\phi$.\med

\noi
\textbf{Example} Let $V$ be a complex linear space with inner product
$\li\,\cdot\,|\,\cdot\,\re$. Let $V'$ denote the dual space of $V$, i.e., the
space of all linear forms $l:V\to\C$. For any $v\in V$, we can define a linear
form $\li v|\in V'$ by $\li v|w:=\li v|w\re$. Then $V'$, together with the map
$v\mapsto\li v|$, is a realization of the conjugate of $V$.\med

If $(V,\pi)$ is a representation of a Lie algebra $\gk$, then we can equip the
conjugate space $\ov V$ with the structure of a representation of the
conjugate Lie algebra $\ov\gk$ by putting
\[
\ov\xb\,\ov v:=\ov{\xb v}.
\]
It is easy to see that this defines a left action of $\ov\gk$ on $\ov V$.
We call $\ov V$, equipped with this left action of $\ov\gk$, the
\emph{conjugate} of the representation $V$.

There is a close relation between Lie algebras and Lie groups. Roughly
speaking, a Lie group is a smooth differentiable manifold that is equipped
with a group structure. In particular, a matrix Lie group $G$ is a group whose
elements are invertible linear operators acting on some finite dimensional
linear space $V$. The Lie algebra of $G$ is then defined as
\[
\hk:=\{A\in\Li(V):e^{tA}\in G\ \forall t\geq 0\}.
\]
In general, this is a real Lie algebra. More generally, one can associate a
Lie algebra to each Lie group (not necessarily a matrix Lie group) and prove
that each Lie algebra is the Lie algebra of some Lie group. Under a certain
condition (simple connectedness), the Lie algebra determines its associated
Lie group uniquely. A finite dimensional representation of a Lie group $G$ is
a pair $(V,\Pi)$ where $V$ is a finite dimensional linear space and
$\Pi:G\to\Li(V)$ is a group homomorphism. Each representation $(V,\pi)$ of a
real Lie algebra $\hk$ gives rise to a representation $(V,\Pi)$ of the
associated Lie group such that $\Pi(e^{tA})=e^{t\pi(A)}$. If $\gk$ is the
complexification of $\hk$ and $(V,\pi)$ is a unitary representation of $\gk$,
then $(V,\Pi)$ is a unitary representation of $G$ in the sense that $\Pi(A)$
is a unitary operator for each $A\in G$. All his is
explained in more detail in Appendix~\ref{A:Lie}.

\subsection{The Lie algebra SU(2)}\label{S:SU2}

The Lie algebra $\suk(2)$ is the three dimensional complex Lie algebra defined
by the commutation relations between its basis elements
\be\label{su2com}
[\sv_{\rm x},\sv_{\rm y}]=2i\sv_{\rm z},\quad
[\sv_{\rm y},\sv_{\rm z}]=2i\sv_{\rm x},\quad
[\sv_{\rm z},\sv_{\rm x}]=2i\sv_{\rm y}.
\ee
It is customary to equip $\suk(2)$ with an adjoint operation that is defined
by
\be\label{su2adj}
\sv_{\rm x}^\ast=\sv_{\rm x},\quad
\sv_{\rm y}^\ast=\sv_{\rm y},\quad
\sv_{\rm z}^\ast=\sv_{\rm z}.
\ee
A faithful unitary representation of $\suk(2)$ is defined by the \emph{Pauli
  matrices}
\be\label{Pauli1}
S_{\rm x}:=\left(\ba{cc}0&1\\1&0\ea\right),\quad
S_{\rm y}:=\left(\ba{cc}0&-i\\i&0\ea\right),\quand
S_{\rm z}:=\left(\ba{cc}1&0\\0&-1\ea\right).
\ee
It is straightforward to check that these matrices are linearly independent
and satisfy the commutation and adjoint relations (\ref{su2com}) and
(\ref{su2adj}). In particular, this shows that $\suk(2)$ is
well-defined.\footnote{Not every set of commutation relations that one can
  write down defines a bona fide Lie algebra. By linearity and skew symmetry,
  specifying $[\xb_i,\xb_j]$ for all $i<j$ uniquely defines a bilinear map
  $[\,\cdot\,,\,\cdot\,]$, but such a map may fail to satisfy the Jacobi
  identity. Similarly, it is not a priori clear that (\ref{su2adj}) defines a
  bona fide adjoint operation, but the faithful unitary representation defined
  by the Pauli matrices shows that it does.}

In general, if $S_{\rm x},S_{\rm y},S_{\rm z}$ are linear operators on
some complex linear space $V$ that satisfy the commutation relations
(\ref{su2com}), and hence define a representation $(V,\pi)$ of $\suk(2)$,
then the so-called \emph{Casimir operator} is defined as
\[
C:=S_{\rm x}^2+S_{\rm y}^2+S_{\rm z}^2.
\]
The operator $C$ is in general not an element of $\{\pi(\xb):\xb\in\suk(2)\}$,
i.e., $C$ does not correspond to an element of the Lie algebra $\suk(2)$. It
does correspond, however, to an element of the so-called \emph{universal
  enveloping algebra} of $\suk(2)$; see Appendix~\ref{A:algLie} below.

The finite-dimensional irreducible representations of $\suk(2)$ are well
understood. Part~(a) of the following proposition follows from
Theorem~\ref{T:Liecomp} in the appendix, using the compactness of the Lie
group ${\rm SU}(2)$. Parts~(b) and (c), and also Proposition~\ref{P:raislow}
below, follow from \cite[Thm~4.32]{Hal03} and a calculation of the Casimir
operator for the representation in Proposition~\ref{P:raislow}.

\bp[Irreducible representations of $\suk(2)$]
Let\label{P:su2clas} $S_{\rm x},S_{\rm y},S_{\rm z}$ be linear operators on a
finite dimensional complex linear space $V$, that satisfy the commutation
relations (\ref{su2com}) and hence define a representation $(V,\pi)$ of
$\suk(2)$. Then:
\begin{itemize}
\item[{\rm\textbf{(a)}}] There exists an inner product
  $\li\,\cdot\,|\,\cdot\,\re$ on $V$, which is unique up to a multiplicative
  constant, such that with respect to this inner product the representation
  $(V,\pi)$ is unitary.
\item[{\rm\textbf{(b)}}] If the representation $(V,\pi)$ is irreducible, then
  there exists an integer $n\geq 1$, which we call the \emph{index} of
  $(V,\pi)$, such that the Casimir operator $C$ is given by $C=n(n+2)I$.
\item[{\rm\textbf{(c)}}] Two irreducible representations $V,W$ of $\suk(2)$
  are equivalent if and only if they have the same index.
\end{itemize}
\ep

Proposition~\ref{P:su2clas} says that the finite dimensional irreducible
representations of $\suk(2)$, up to isomorphism, can be labeled by their index
$n$, which is a natural number $n\geq 1$. We next describe what an irreducible
representation with index $n$ looks like. In spite of the beautiful symmetry
of the commutation relations (\ref{su2com}), it will be useful to work with a
different, less symmetric basis $\{\jb^-,\jb^+,\jb^0\}$ defined as
\be\label{jdef}
\jb^-:=\ha(\sv_{\rm x}-i\sv_{\rm y}),\quad
\jb^+:=\ha(\sv_{\rm x}+i\sv_{\rm y}),\quand
\jb^0:=\ha\sv_{\rm z},
\ee
which satisfies the commutation and adjoint relations:
\be\label{railow}
[\jb^0,\jb^\pm]=\pm \jb^\pm,\quad[\jb^-,\jb^+]=-2\jb^0,\quad
(\jb^-)^\ast=\jb^+,\quad(\jb^0)^\ast=\jb^0.
\ee
The next proposition describes what an irreducible representation of $\suk(2)$
with index $n$ looks like.

\bp[Raising and lowering operators]
Let\label{P:raislow} $V$ be a finite dimensional complex linear space that is
equipped with an inner product and let $J^\pm,J^0$ be linear operators on $V$
that satisfy the commutation and adjoint relations (\ref{railow}) and hence
define a unitary representation $(V,\pi)$ of $\suk(2)$. Assume that $(V,\pi)$
is irreducible and has index $n$. Then $V$ has dimension $n+1$ and there exists
an orthonormal basis
\[
\{\phi(-n/2),\phi(-n/2+1),\ldots,\phi(n/2)\}
\]
such that
\bc\label{raislow}
\dis J^0\phi(k)&=&\dis k\phi(k),\\[5pt]
\dis J^-\phi(k)&=&\dis\sqrt{(n/2-k+1)(n/2+k)}\phi(k-1),\\[5pt]
\dis J^+\phi(k)&=&\dis\sqrt{(n/2-k)(n/2+k+1)}\phi(k+1)
\ec
for $k=-n/2,-n/2+1,\ldots,n/2$, with the conventions $J^-\phi(-n/2):=0$
and $J^+\phi(n/2):=0$.
\ep

We see from (\ref{raislow}) that $\phi(k)$ is an eigenvector of $J^0$ with
eigenvalue $k$, and that the operators $J^\pm$ maps such an eigenvector into
an eigenvector with eigenvalue $k\pm1$, respectively. In view of this, $J^\pm$
are called \emph{raising} and \emph{lowering} operators, or also
\emph{creation} and \emph{annihilation} operators. It is
instructive to see how this property of $J^\pm$ follows rather easily from the
commutation relations (\ref{railow}). Indeed, if $\phi(k)$ is an eigenvector
of $J^0$ with eigenvalue $k$, then the commutation relations imply that
\[
J^0J^+\phi(k)=\big(J^+J^0+[J^0,J^+]\big)\phi(k)
=\big(J^+J^0+J^+\big)\phi(k)=(k+1)J^+\phi(k),
\]
which shows that $J^+\phi(k)$ is a (possibly zero) multiple of $\phi(k+1)$.
The concept of raising and lowering operators can be generalized to other Lie
algebras.

\subsection{The Lie algebra SU(1,1)}\label{S:SU11}

The Lie algebra $\suk(1,1)$ is defined by the commutation relations
\be\label{su11com}
[\tb_{\rm x},\tb_{\rm y}]=2i\tb_{\rm z},\quad
[\tb_{\rm y},\tb_{\rm z}]=-2i\tb_{\rm x},\quad
[\tb_{\rm z},\tb_{\rm x}]=2i\tb_{\rm y}.
\ee
Note that this is the same as (\ref{su2com}) except for the minus sign in the
second equality. A faithful representation is defined by the matrices
\be\label{psPauli}
T_{\rm x}:=\left(\ba{cc}0&1\\-1&0\ea\right),\quad
T_{\rm y}:=\left(\ba{cc}0&i\\i&0\ea\right),\quad
T_{\rm z}:=\left(\ba{cc}1&0\\0&-1\ea\right).
\ee
It is customary to equip $\suk(1,1)$ with an adjoint operation such that
\be\label{su11adj}
\tb_{\rm x}^\ast=\tb_{\rm x},\quad
\tb_{\rm y}^\ast=\tb_{\rm y},\quad
\tb_{\rm z}^\ast=\tb_{\rm z}.
\ee
Note however, that the matrices in (\ref{psPauli}) are not self-adjoint and
hence do not define a unitary representation of $\suk(1,1)$. In fact, all
unitary irreducible representations of $\suk(1,1)$ are infinite
dimensional.
\footnote{Since $\suk(1,1)$ is simple, all representations are
    faithful. As explained in Subsection~\ref{S:gralg}, each Lie algebra is
    the Lie algebra of a unique simply connected Lie group. In the case of
    $\suk(1,1)$, this is the universal cover $H$ of the Lie group ${\rm
      SU}(1,1)$ (the latter itself not being simply connected). By
    Theorem~\ref{T:simpcon} in the appendix, each representation of
    $\suk(1,1)$ gives rise to a representation of $H$. Since $H$ is not
    compact, the existence of a finite dimensional unitary representation
    would now contradict Lemma~\ref{L:noncomp} in the appendix.} In a given
representation of $\suk(1,1)$, the \emph{Casimir operator} is defined as
\be\label{CasimT}
C:=(\ha T_{\rm x})^2-(\ha T_{\rm y})^2-(\ha T_{\rm z})^2.
\ee
Again, it is useful to introduce raising and lowering operators, defined as
\[
\kb^0:=\ha\tb_{\rm x}\quand\kb^\pm:=\ha(\tb_{\rm y}\pm i\tb_{\rm z}),
\]
which satisfy the commutation and adjoint relations
\be\label{kpm}
[\kb^0,\kb^\pm]=\pm\kb^\pm,\quad[\kb^-,\kb^+]=2\kb^0,\quad
(\kb^-)^\ast=\kb^+,\quad(\kb^0)^\ast=\kb^0,
\ee
The following proposition is rewritten from \cite[formulas (8) and
  (9)]{Nov04}, where this is stated without proof or reference.
The constant $r>0$ below is called the \emph{Bargmann index} \cite{Bar47,Bar61}.

\bp[Representations of $\suk(1,1)$]
For\label{P:su11rep} each real constant $r>0$, there exists an irreducible
unitary representation of $\suk(1,1)$ on a Hilbert space with orthonormal basis
$\{\phi(0),\phi(1),\ldots\}$ on which the operators $K^0,K^\pm$ act as
\bc\label{Krep}
\dis K^0\phi(k)&=&\dis(k+r)\phi(k),\\[5pt]
\dis K^-\phi(k)&=&\dis 1_{\{k\geq 1\}}\sqrt{k(k-1+2r)}\phi(k-1),\\[5pt]
\dis K^+\phi(k)&=&\dis\sqrt{(k+1)(k+2r)}\phi(k+1).
\ec
In this representation, the Casimir operator is given by $C=r(r-1)I$.
\ep

In what follows, we will need one more representation of $\suk(1,1)$, as well as
a representation of its conjugate Lie algebra. Fix $\al>0$ and consider the
following operators acting on smooth functions $f:\half\to\R$:
\be\label{Kidef}\bac
\dis\Ki^-f(z)&=&\dis z\diff{z}f(z)+\al\dif{z}f(z),\\[5pt]
\dis\Ki^+f(z)&=&\dis zf(z),\\[5pt]
\dis\Ki^0f(z)&=&\dis z\dif{z}f(z)+\ha\al f(z).
\ec
One can check that these operators satisfy the commutation relations
(\ref{kpm}) of the Lie algebra $\suk(1,1)_\C$, i.e.,
\be\label{Kcomrel}
[\Ki^0,\Ki^\pm]=\pm\Ki^\pm\quand[\Ki^-,\Ki^+]=2\Ki^0,
\ee
and hence define a representation of $\suk(1,1)$. One can check that the
  Casimir operator \eqref{CasimT} for this representation is
  $C=\frac{\alpha}2(\frac\alpha2-1)I$ and hence the Bargmann index is
  $r=\al/2$.

Next, fix again $\al>0$ and consider the following operators acting on functions
$f:\N\to\R$:
\be\bac\label{Kdef}
\dis K^-f(x)&=&\dis xf(x-1),\\[5pt]
\dis K^+f(x)&=&\dis (\al+x)f(x+1),\\[5pt]
\dis K^0f(x)&=&\dis (\ha\al+x)f(x).
\ec
One can check that these operators satisfy the commutation relations
\be\label{ducomrel}
[K^\pm,K^0]=\pm K^\pm\quand[K^+,K^-]=2K^0.
\ee
This is similar to (\ref{kpm}), except that the order of the elements inside
the commutator is reversed. In view of the remarks at the end of
Section~\ref{S:Lieint}, this means that the operators $K^0,K^\pm$ define a
representation of the conjugate Lie algebra associated with $\suk(1,1)$.
We will see in Section~\ref{S:SIP} below that the conjugate of the
representation in (\ref{Kdef}) is equivalent to the representation in
(\ref{Kidef}), provided we choose for both the same $\al$.

A complete classification of all irreducible representations of $\suk(1,1)$,
including infinite dimensional ones, is described in the book
\cite{VK91}.\footnote{The monumental encyclopedic book \cite{VK91} is written
  in a style that some readers may need to get used to, since it does not use
  the usual theorem-proof layout but rather states an enormous amount of facts
  in the main text while leaving a lot of detail to be filled in by the
  reader.}

\subsection{The Heisenberg algebra}\label{S:Heis}

The \emph{Heisenberg algebra} $\hk$ is the three dimensional complex Lie algebra
defined by the commutation relations
\be\label{Heiscom}
[\ab^-,\ab^+]=\ab^0,\quad[\ab^-,\ab^0]=0,\quad[\ab^+,\ab^0]=0.
\ee
It is customary to equip $\hk$ with an adjoint operation that is defined by
\be
(\ab^\pm)^\ast=\pm\ab^\pm,\quad(\ab^0)^\ast=\ab^0.
\ee
The \emph{Schr\"odinger representation} of $\hk$ is defined by
\be\label{Schroed}
A^-f(x)=\dif{x}f(x),\quad A^+f(x)=xf(x),\quad A^0f(x)=f(x),
\ee
which are interpreted as operators on the Hilbert space $L^2(\R,\di x)$ of
complex functions on $\R$ that are square integrable with respect to the
Lebesgue measure. Note in this representation, $A^0$ is the identity
operator. Any representation of $\hk$ with this property is called a
\emph{central} representation.\footnote{More generally, the \emph{center} of a
  Lie algebra $\gk$ is the linear space $\ck:=\{\cb\in\gk:[\xb,\cb]=0\ \forall
  \xb\in\gk\}$. A \emph{central} representation of a Lie algebra is then a
  representation $(V,\pi)$ such that for each $\cb\in\ck$, there exists a
  $c\in\C$ such that $\pi(\cb)=cI$. Note that with this definition, if
  $(V,\pi)$ is a faithful central representation of $\hk$, then we can always
  ``normalize'' it by multiplying $\pi$ with a constant so that
  $\pi(\ab^0)=I$.} The Schr\"odinger representation is a
unitary representation, i.e., $A^-$ is skew symmetric and $A^+$ and $A^0$ are
self-adjoint, viewed as linear operators on the Hilbert space $L^2(\R,\di x)$.

Since $iA^-$ and $A^+$ are self-adjoint, by Stone's theorem, 
one can define collections of unitary operators $(U^-_t)_{t\in\R}$ and
$(U^+_t)_{t\in\R}$ by
\be
U^-_s:=\ex{tA^-}\quand U^+_t:=\ex{itA^+}.
\ee
These operators form one-parameter groups in the sense that $U^\pm_0=I$ and
$U^\pm_sU^\pm_t=U^\pm_{s+t}$ $(s,t\in\R)$.  Note that we have a factor $i$ in
the definition of $U^+_t$ but not in the definition of $U^-_s$, because $A^+$
is self-adjoint but $A^-$ is skew symmetric. The commutation relations
(\ref{Heiscom}) lead, at least formally, to the following commutation relation
between $U^-_s$ and $U^+_t$
\be\label{SNcom}
U^-_sU^+_t=e^{ist}U^+_tU^-_s\qquad(s,t\in\R).
\ee
Indeed, for small $\eps$, we have
\be\ba{l}
\dis U^-_{\eps s}U^+_{\eps t}\\[5pt]
\dis\ =\big(I+\eps sA^-+\ha\eps^2 s^2(A^-)^2+O(\eps^3)\big)
\big(I+i\eps tA^+-\ha\eps^2 t^2(A^+)^2+O(\eps^3)\big)\\[5pt]
\dis\ =I+\eps sA^-+\ha\eps^2 s^2(A^-)^2+i\eps tA^+-\ha\eps^2 t^2(A^+)^2
+i\eps^2 stA^-A^++O(\eps^3)\\[5pt]
\dis\ =I+\eps sA^-+\ha\eps^2 s^2(A^-)^2+i\eps tA^+-\ha\eps^2 t^2(A^+)^2
+i\eps^2 stA^+A^-\\[5pt]
\dis\hspace{8cm}
+i\eps^2 st[A^-,A^+]+O(\eps^3)\\[0pt]
\dis\ =\big(1+i\eps^2 st+O(\eps^3)\big)U^+_{\eps t}U^-_{\eps s}+O(\eps^3).
\ec
The commutation relation (\ref{SNcom}) then follows formally by writing
\bc
\dis U^-_sU^+_t
&=&\dis(U^-_{s/n})^n(U^+_{t/n})^n\\[5pt]
&=&\dis\big(1+in^{-2}st+O(n^{-3})\big)^{n^2}(U^+_{t/n})^n(U^-_{s/n})^n
\asto{n}e^{ist}U^+_tU^-_s.
\ec
The \emph{Stone-von Neumann theorem} states that all unitary, central
representations of the Heisenberg algebra that satisfy (\ref{SNcom}) are
equivalent \cite{Ros04}.

\subsection{The direct sum and the tensor product}\label{S:sumprod}

If $V$ is a linear space and $V_1,\ldots,V_n$ are linear
subspaces of $V$ such that every element $v\in V$ can uniquely be written as
\[
v=v_1+\cdots+v_n
\]
with $v_i\in V_i$, then we say that $V$ is the \emph{direct sum} of
$V_1,\ldots,V_n$ and write $V=V_1\oplus\cdots\oplus V_n$.
If $\Om_1,\Om_2$ are finite sets and $\C^{\Om_1}$ denotes the linear space of
all functions $f:\Om_i\to\C$, then we have the natural isomorphism
\[
\C^{\Om_1\uplus\Om_2}\cong\C^{\Om_1}\oplus\C^{\Om_2},
\]
where $\Om_1\uplus\Om_2$ denotes the disjoint union of $\Om_1$ and $\Om_2$.

If $\gk_1,\ldots,\gk_n$ are Lie algebras, then we equip
$\gk_1\oplus\cdots\oplus\gk_n$ with the structure of a Lie algebra by putting,
for $\xb_i,\yb_i\in\gk_i$,
\be\label{Liesum}
\big[\xb_1+\cdots+\xb_n,\yb_1+\cdots+\yb_n\big]
:=[\xb_i,\yb_i]+\cdots+[\xb_n,\yb_n].
\ee
Note that this has the effect that elements of diffent Lie algebras
$\gk_1,\ldots,\gk_n$ mutually commute. In particular, if
$\{\xb^1_1,\xb^2_1,\xb^3_1\}$ and $\{\xb^1_2,\xb^2_2,\xb^3_2\}$ are bases for
$\gk_1$ and $\gk_2$, respectively, then
\[
\{\xb^1_1,\xb^2_1,\xb^3_1,\xb^1_2,\xb^2_2,\xb^3_2\}
\]
is a basis for $\gk_1\oplus\gk_2$ and $[\xb^k_i,\xb^l_j]=0$ whenever $i\neq j$.

By definition, a \emph{bilinear} map of two variables is a function that is
linear in each of its arguments. If $V$ and $W$ are finite dimensional linear
spaces, then their \emph{tensor product} is a linear space $V\otimes W$
together with a bilinear map
\[
V\times W\ni(v,w)\mapsto v\otimes w\in V\otimes W
\]
that has the property:
\begin{quote}
If $F$ is another linear space and $b:V\times W\to
F$ is bilinear, then there exists a
unique linear map $\ov b:V\otimes W\to F$ such that
\[
\ov b(v\otimes w)=b(v,w)\qquad(v\in V,\ w\in W).
\]
\end{quote}
The tensor product of three or more spaces is defined similarly.
One can show that all realizations of the tensor product are naturally
isomorphic. If $\{e(1),\ldots,e(n)\}$ and $\{f(1),\ldots,f(m)\}$ are bases for
$V$ and $W$, then one can prove that
\be\label{etimf}
\big\{e(i)\otimes f(j):1\leq i\leq n,\ 1\leq j\leq m\big\}
\ee
is a basis for $V\otimes W$. In particular, this means that one has the
natural isomorphism
\be\label{tensfunc}
\C^{\Om_1\times\Om_2}\cong\C^{\Om_1}\otimes\C^{\Om_2}.
\ee
If $A\in\Li(V)$ and $B\in\Li(V)$, then one defines $A\otimes B\in\Li(V\otimes
W)$ by
\be\label{AotB}
(A\otimes B)(v\otimes w):=(Av)\otimes(Bw).
\ee
We note that not every element of $V\otimes W$ is of the form
$v\otimes w$ for some $v\in V$ and $w\in W$. Nevertheless, since the
right-hand side of (\ref{AotB}) is bilinear in $v$ and $w$, the defining 
property of the tensor product tells us that this formula unambiguously
defines a linear operator on $V\otimes W$.

One can check that the notation $A\otimes B$ is good notation in the sense
that the space $\Li(V\otimes W)$ together with the bilinear map $(A,B)\mapsto
A\otimes B$ is a realization of the tensor product $\Li(V)\otimes\Li(W)$.
Thus, one has the natural isomorphism
\[
\Li(V\otimes W)\cong\Li(V)\otimes\Li(W).
\]
If $V$ and $W$ are equipped with inner products,
then we equip $V\otimes W$ with an inner product by putting
\be\label{tensin}
\li v\otimes w|\eta\otimes\xi\re
:=\li v|\eta\re\li w|\xi\re,
\ee
which has the effect that if $\{e(1),\ldots,e(n)\}$ and $\{f(1),\ldots,f(m)\}$
are orthonormal bases for $V$ and $W$, then the basis for $V\otimes W$ in
(\ref{etimf}) is also orthonormal. Again, one needs the defining 
property of the tensor product to see that (\ref{tensin}) is a good definition.

If $V,W$ are representations of Lie algebras $\gk,\hk$, respectively, then we
can naturally equip the tensor product $V\otimes W$ with the structure of a
representation of $\gk\oplus\hk$ by putting
\be\label{gplush}
(\xb+\yb)(v\otimes w):=(\xb v)\otimes(\yb w).
\ee
Again, since the right-hand side is bilinear, using the defining property of
the tensor product, one can see that this is a good definition.

Let $V_1,V_2$ be representations of some Lie algebra $\gk$,
and let $W_1,W_2$ be representations of another Lie algebra $\hk$.
Let $\phi:V_1\to V_2$ and $\psi:W_1\to W_2$ be intertwiners. Then one can
check that
\be\label{prodtwine}
\phi\otimes\psi:V_1\otimes W_1\to V_2\otimes W_2
\ee
is also an intertwiner.

If $\hk_1,\ldots,\hk_n$ are $n$ copies of the Heisenberg algebra, and
$\ab^-_i,\ab^+_i,\ab^0_i$ are basis elements of $\hk_i$ that satisfy the
commutation relations (\ref{Heiscom}), then a basis for
$\hk_1\oplus\cdots\oplus\hk_n$ is formed by all elements $\ab^\pm_i,\ab^0_i$
with $i=1,\ldots,n$, and these satisfy
\[
[\ab^-_i\ab^+_j]=\de_{ij}\ab^0_i\quand[\ab^\pm_i,\ab^0_j]=0.
\]
Since the center of $\hk_1\oplus\cdots\oplus\hk_n$ is spanned by the elements
$\ab^0_i$ with $i=1,\ldots,n$, a central representation of
$\hk_1\oplus\cdots\oplus\hk_n$ must map all these elements to multiples of the
identity. In particular, a central representation of
$\hk_1\oplus\cdots\oplus\hk_n$ is never faithful (unless $n=1$).
The Lie algebra $\hk(n)$ is the $2n+1$ dimensional Lie algebra with basis
elements $\ab^\pm_i$ $(i=1,\ldots,n)$ and $\ab^0$, which satisfy the
commutation relations
\[
[\ab^-_i\ab^+_j]=\de_{ij}\ab^0\quand[\ab^\pm_i,\ab^0]=0.
\]
A \emph{central} representation of $\hk(n)$ is a representation $(V,\pi)$
such that $\pi(\ab^0)=I$. The Schr\"odinger representation of the
``$n$-dimensional'' Heisenberg algebra is the central representation of
$\hk(n)$ on $L^2(\R^n,\di x)$ given by
\be
A^-f(x)=\dif{x_i}f(x)\quand A^+f(x):=x_if(x).
\ee

\section{The algebraic approach to duality}\label{S:repdu}

After our excursion into the theory of Lie algebras, we return to our main
topic, which is the algebraic approach to Markov process duality. We recall
from Lemma~\ref{L:build} that if a Markov generator $L$ can be written in
terms of ``building blocks'' $A_i$ that each have a dual $B_i$ with respect to
some duality function $D$, then also $L$ has a dual $\hat L$ with respect to
$D$. As mentioned at the end of Section~\ref{S:alg}, it may be a good idea to
choose the $A_i$'s so that they define a representation of some Lie algebra.
The next proposition says that in such a situation, other, equivalent
representations of the same Lie algebra may lead to dual Markov processes.

Recall the definition of a conjugate Lie algebra $\ov\gk$ from
Section~\ref{S:Lieint}. If $Y_1,\ldots,Y_n$ are matrices that define
a representation of $\ov\gk$, then their adjoints $Y_1^\dgg,\ldots,Y_n^\dgg$
define a representation of the original Lie algebra $\gk$.

\bp[Intertwiners as duality functions]
Let\label{P:genprin} $L$ be the generator of a Markov process with finite
state space $\Om$. Let $X_1,\ldots,X_n$ be linear operators on $\C^\Om$ that
form a representation of some Lie algebra $\gk$. Assume that $L$
can be written as a linear combination of products of the operators
$X_1,\ldots,X_n$
\be
L=\sum_{(i_1,\ldots,i_k)\in\Ii}r_{i_1,\ldots,i_k}X_{i_1}\cdots X_{i_k},
\ee
where $\Ii$ is some finite set whose elements are sequences
$(i_1,\ldots,i_k)$ with $k\geq 0$ and $1\leq i_m\leq n$ for each $m$.
Assume that $Y_1,\ldots,Y_n$ are linear operators on $\C^{\hat\Om}$ that define
a representation of the conjugate Lie algebra $\ov\gk$. Assume that the
representations of $\gk$ defined by $Y_1^\dgg,\ldots,Y_n^\dgg$ and
$X_1,\ldots,X_n$ are equivalent. Then there is a bijective intertwiner $D$, 
i.e., $X_iD=DY^\dgg_i$ for each $i$, and $L$ is dual w.r.t.\ the duality
function $D$ to the operator
\be\label{LXX}
\hat L:=\sum_{(i_1,\ldots,i_k)\in\Ii}r_{i_1,\ldots,i_k}Y_{i_k}\cdots Y_{i_1}.
\ee
\ep
\bpro
By definition, two representations are equivalent if and only if there exists
a bijective intertwiner. The fact that $L$ is dual to the operator in
(\ref{LXX}) is then immediate from Lemma~\ref{L:build}.
\epro

At first sight, it may seem unlikely that Proposition~\ref{P:genprin} could be
of much use. Even if we can write a generator in terms of a basis of a
representation of some Lie algebra $\gk$, and we also find some representation
of the conjugate Lie algebra $\ov\gk$, we still have to be lucky in the sense
that the representations of $\gk$ defined by $Y_1^\dgg,\ldots,Y_n^\dgg$ and
$X_1,\ldots,X_n$ are equivalent, and there is no guarantee that the operator
in (\ref{LXX}) is a Markov generator. Nevertheless, in what follows, we will
see that Proposition~\ref{P:genprin} can help us find nontrivial dualities.
In the next subsection, we demonstrate this on the operator $L$ from
(\ref{WF}), which is the generator of a Wright-Fisher diffusion with
selection.

\subsection{Self-duality of the Wright-Fisher diffusion}\label{S:WFself}

In Subsection~\ref{S:alg}, we have seen that the operator $L$ from (\ref{WF})
can as in (\ref{LAexpr}) be written in terms of the ``building blocks''
$A^\pm$ from (\ref{Apm}). As we have seen in (\ref{ApmHeis}), these operators
satisfy
\be\label{AHeis}
[A^-,A^+]=I,
\ee
and hence define a central representation of the Heisenberg algebra $\hk$, as
defined in Subsection~\ref{S:Heis}.

It will be convenient to find a way of writing $L$ in a more symmetric way
than in (\ref{LAexpr}). To this aim, we change the definitions of $A^\pm$ to
\be\label{Apm2}
A^-f(x):=\frac{-1}{\sqrt{s}}\dif{x}f(x)\quand A^+f(x):=\sqrt{s}xf(x),
\ee
which again satisfy (\ref{AHeis}), and we write $L$ in terms of these new
building blocks as
\be\label{LAexpr2}
L=-A^+(\sqrt{s}-A^+)A^-(\sqrt{s}-A^-).
\ee
We observe from (\ref{AHeis}) that setting $B^-:=A^+$ and $B^+:=A^-$ defines
operators such that $[B^-,B^+]=-I$, i.e., $B^-,B^+$ define a central
representation of the conjugate Heisenberg algebra $\ov\hk$.

We recall from Section~\ref{S:Heis} that the Stone-von Neumann theorem states
that, more or less, all central representations of the Heisenberg algebra are
equivalent. In view of this and Proposition~\ref{P:genprin}, we may expect
that the operator
\be\label{hatLAexpr}
\hat L=-(\sqrt{s}-B^-)B^-(\sqrt{s}-B^+)B^+
\ee
is dual to $L$ with respect to some (so far unknown) duality function $D$.
(Here (\ref{hatLAexpr}) is obtained from (\ref{LAexpr2}) by replacing $A^\pm$
by $B^\pm$ and reversing the order of the factors.) Since $B^\pm=A^\mp$, we
observe that in fact $\hat L=L$, so our calculations lead us to suspect that the
Wright-Fisher diffusion with selection parameter $s>0$ should be self-dual.

We still need to find the duality function $D$. This function must satisfy
\be
\frac{-1}{\sqrt{s}}\dif{x}D(x,y)=A^-D(\,\cdot\,,y)(x)
=B^-D(x,\,\cdot\,)(y)=\sqrt{s}yD(x,y),
\ee
which says that $\dif{x}D(x,y)=-syD(x,y)$ and leads to the requirement that
$D(x,y)=D(0,y)e^{-syx}$. In a similar way, the requirement $A^+D=DB^+$ yields
$D(x,y)=D(x,0)e^{-syx}$ and in particular $D(0,y)=D(0,0)$. Thus, we find that
up to a multiplicative constant, there is a unique duality function, which is
given by
\be\label{Dexp}
D(x,y)=\ex{-sxy}\qquad\big(x,y\in[0,1]\big),
\ee
and we conclude that the Wright-Fisher diffusion with selection parameter
$s>0$ is self-dual with this duality function.

The argument above was heuristic, but quite smooth. What is remarkable about
it is that while usually, the discovery of a duality starts with a clever
choice for the duality function, here, the duality function came at the very
end. Hidden behind this is the Stone-von Neumann theorem which says that two
``good'' representations of the Heisenberg algebra must necessarily be
equivalent. We did not check the conditions of this theorem in detail (this is
why the argument is only heuristic), but rather used it as an inspiration.  A
priori, there was no guarantee that the operator in (\ref{hatLAexpr}) would be
a Markov generator, but since $\hat L=L$ and $L$ is a Markov generator, this
turned out right as well.

\noi
\textbf{Remark} It is possible to ``discover'' the moment dual (\ref{WFmomdu})
of the Wright-Fisher duality along similar lines as we have discovered its
self-duality here, by considering a suitable representation of the conjugate
Heisenberg algebra $\ov\hk$ on functions $f:\N\to\R$ and applying
Propositions~\ref{P:genprin} and \ref{P:Hilb}. Such a derivation is less
natural, however, since it requires choosing a rather peculiar representation
of $\ov\hk$ that more or less has the duality function from (\ref{Dmom})
tacitly built into it.

\subsection{Intertwiners and duality functions}\label{S:contwine}

In the previous subsection, just before (\ref{hatLAexpr}) we appealed to
Proposition~\ref{P:genprin}. In doing so, we cheated in the sense that the
operators $A^\pm$ from (\ref{Apm2}) do not act on a finite-dimensional
space. The most obvious consequence of this is that it is not clear how the
adjoint operators $B_i^\dgg$ from Proposition~\ref{P:genprin} should be
defined. Closely related to this is that in the infinite dimensional setting,
it is not immediately clear that duality functions define intertwiners and
vice versa. In this subsection we show that these difficulties can be resolved
by introducing a suitable inner product on the spaces of complex functions on
$\Om$ and $\hat\Om$, respectively.

Assume that $X_1,\ldots,X_n$ and $Y_1,\ldots,Y_n$ are
linear operators on $L^2$-spaces $L^2(\Om,\mu)$ and $L^2(\hat\Om,\nu)$,
respectively, that define representations of a Lie algebra $\gk$ and its
conjugate $\ov\gk$. Let $Y^\ast_i$ denote the adjoint of $Y_i$ with
  respect to the inner product on $L^2(\hat\Om,\nu)$. Assume that
$\Phi:L^2(\hat\Om,\nu)\to L^2(\Om,\mu)$ is a linear operator of the form
\be\label{Phikern}
\Phi g(x)=\int g(y)D(x,y)\nu(\di y),
\ee
for some function $D:\Om\times\hat\Om\to\C$ such that the expressions in
(\ref{XYD}) below are well-defined.

\bp[Intertwiners and duality functions]
The\label{P:Hilb} operator $\Phi$ is an intertwiner of the
representations defined by $X_1,\ldots,X_n$ and
$Y^\ast_1,\ldots,Y^\ast_n$, i.e.,
\[
X_i\Phi=\Phi Y^\ast_i\quad(i=1,\ldots,n),
\]
if and only if $D$ is a duality function, in the sense that
\be\label{XYD}
X_iD(\,\cdot\,,y)(x)=Y_iD(x,\,\cdot\,)(y)\qquad(i=1,\ldots,n)
\ee
for a.e.\ $x,y$ with respect to the product measure $\mu\otimes\nu$.
\ep
\bpro
We observe that
\[\ba{l}
\dis\int\ov{f(x)}\mu(\di x)\int g(y)\nu(\di y)X_iD(\,\cdot\,,y)(x)
=\int g(y)\nu(\di y)\li f|X_iD(\,\cdot\,,y)\re_\mu\\[5pt]
\dis\quad=\int g(y)\nu(\di y)\li X^\ast_if|D(\,\cdot\,,y)\re_\mu
=\int\ov{X^\ast_if(x)}\mu(\di x)\int g(y)\nu(\di y)D(x,y)\\[5pt]
\dis\quad=\li X^\ast_if|\Phi g\re_\mu=\li f|X_i\Phi g\re_\mu
\ea\]
and
\[\ba{l}
\dis\int\ov{f(x)}\mu(\di x)\int g(y)\nu(\di y)Y_iD(x,\,\cdot\,)(y)
=\int\ov{f(x)}\mu(\di x)\li g|Y_iD(x,\,\cdot\,)\re_\nu\\[5pt]
\dis\quad=\int\ov{f(x)}\mu(\di x)\li Y^\ast_ig|D(x,\,\cdot\,)\re_\nu
=\int\ov{f(x)}\mu(\di x)\int Y^\ast_ig(y)\nu(\di y)D(x,y)\\[5pt]
\dis\quad=\li f|\Phi Y^\ast_ig\re_\mu.
\ea\]
Since this holds for all $f,g$, the statement follows.
\epro

\noi
\textbf{Remark} Proposition \ref{P:Hilb} allows us to obtain an intertwiner
from a duality function. Conversely, if $\Phi:L^2(\hat\Om,\nu)\to
L^2(\Om,\mu)$ is a bounded linear operator, then setting
\[
\Delta(f\otimes g) :=\int f(x)\Phi g(x)\,\mu(\di x)
\]
defines a linear form on the linear span of all functions of the form
$f\otimes g$. If $\De$ is bounded,\footnote{Using Cauchy-Schwarz,
it is easy to see that $|\De(f\otimes g)|\leq\|\Phi\|\,\|f\otimes g\|$,
proving that $\De$ is bounded on functions of the form $f\otimes
g$. Nevertheless, $\De$ may fail to be bounded on the linear span of such
functions. A counterexample is $\Om=\hat\Om=[0,1]$, $\mu=\nu=$ Lebesgue
measure, and $\Phi$ the identity map, which gives $\De(F)=\int_0^1F(x,x)\,\di
x$. Since the Lebesgue measure on the diagonal $\{(x,y):x=y\}$ does not have a
density w.r.t.\ $\mu\otimes\nu$, this does not correspond to a bounded linear
form on $L^2(\Om\times\hat\Om,\mu\otimes\nu)$.} then it can uniquely be
extended to a bounded linear form on
\[
L^2(\Om,\mu)\otimes L^2(\hat\Om,\nu)\cong L^2(\Om\times\hat\Om,\mu\otimes\nu),
\]
so that by the Riesz representation theorem there exists a $D\in
L^2(\Om\times\hat\Om,\mu\otimes\nu)$ such that
\[
\Delta(f\otimes g) :=\int f(x)D(x,y)g(x)\,\mu(\di x)\nu(\di y),
\]
proving that $\Phi$ is of the form \eqref{Phikern} (although
there is no guarantee that $D(\,\cdot\,,y)$ and $D(x,\,\cdot\,)$ are in the
domains of $X_i$ and $Y_i$, resp., if these are unbounded operators).

\subsection{The symmetric exclusion process}\label{S:SEP}

In this subsection, we demonstrate Proposition~\ref{P:genprin} on a simple
example, which involves the simple exclusion process and the Lie algebra
$\suk(2)$. In the end, we find a self-duality that is not entirely trivial,
but also not very useful. The present subsection serves mainly as a warm-up
for Subsection~\ref{S:SIP} where we will replace $\suk(2)$ by $\suk(1,1)$.

Let $S$ be a finite set and let $r:S\times S\to\half$ be a function that is
symmetric in the sense that $r(i,j)=r(j,i)$. Consider the Markov process with
state space $\Om=\{0,1\}^S$ and generator
\be
Lf(x):=\sum_{ij}r(i,j)1_{\{(x_i,x_j)=(1,0)\}}\big\{f(x-\de_i+\de_j)-f(x)\big\},
\ee
where $\de_i\in\Om$ is defined as $\de_i(j):=1_{\{i=j\}}$.
Then $L$ is the generator of a \emph{symmetric exclusion process} or
\emph{SEP}. We define operators $J^\pm_i$ and $J^0_i$ by
\be\ba{l}
\dis J^-_if(x):=1_{\{x_i=0\}}f(x+\de_i),\quad
J^+_if(x):=1_{\{x_i=1\}}f(x-\de_i),\\[5pt]
\quand J^0_if(x):=(x_i-\ha)f(x).
\ec
It is straightforward to check that 
\be\label{excom}
[J^0_i,J^\pm_j]=\pm\de_{ij} J^\pm_i\quand
[J^-_i,J^+_j]=-2\de_{ij}J^0_i.
\ee
It follows that the operators $J^\pm_i$ and $J^0_i$ define a representation of
a Lie algebra that consists of a direct sum of copies of $\suk(2)$, with one
copy for each site $i\in S$. We can write the generator $L$ of the symmetric
exclusion process in terms of the operators $J^\pm_i$ and $J^0_i$ as
\be\label{excl1}
L=\sum_{\{i,j\}}r(i,j)\big[J^-_iJ^+_j+J^-_jJ^+_i+2J^0_iJ^0_j-\ha I\big],
\ee
where we are summing over all unordered pairs $\{i,j\}$. We observe that the
operators
\be\label{KJdef}
K^\pm_i:=J^\pm_i,\quand K^0_i:=-J^0_i
\ee
satisfy the same commutation relations as $J^\pm_i$ and $J^0_i$, except that
each commutation relation gets an extra minus sign. This shows that the
operators $K^\pm_i$ and $K^0_i$ define a representation of the conjugate Lie
algebra $\ov{\suk(2)}$. Moreover, we can alternatively write the generator in
(\ref{excl1}) as
\be\label{excl2}
L=\sum_{\{i,j\}}r(i,j)\big[K^+_jK^-_i+K^+_iK^-_j+2K^0_jK^0_i-\ha I\big].
\ee

We recall from Subsection~\ref{S:SU2} that two irreducible representations of
$\suk(2)$ with the same dimension are necessarily equivalent. In view of this,
we conjecture that there should exist an intertwiner $D$, unique up to a
multiplicative constant, such that $J^\pm_iD=D(K^\pm_i)^\dgg$ and
$J^0_iD=D(K^0_i)^\dgg$ for all $i$. By the general principle in
Proposition~\ref{P:genprin}, such an intertwiner is a self-duality function
for the symmetric exclusion process.

We observe that all our operators act on the space of all complex functions on
$\{0,1\}^S$, which in view of (\ref{tensfunc}) is given by
\be\label{tensC}
\C^{\,\txt\{0,1\}^S}\cong\bigotimes_{i\in S}\C^{\{0,1\}}.
\ee
For example, if $S=\{1,2,3\}$ consists of only three sites, then in line with
(\ref{gplush}),
\[
J^0_1=J^0\otimes I\otimes I,\quad
J^0_2=I\otimes J^0\otimes I,\quand
J^0_3=I\otimes I\otimes J^0,
\]
and similarly for $J^\pm_1,J^\pm_2$, and $J^\pm_3$. Here
\bc\label{Jpmo}
\dis J^-f&=&\dis\left(\ba{cc}0&0\\ 1&0\ea\right)
\left(\ba{c}f(1)\\ f(0)\ea\right)
=\left(\ba{c}0\\ f(1)\ea\right),\\[17pt]
\dis J^+f&=&\dis\left(\ba{cc}0&1\\ 0&0\ea\right)
\left(\ba{c}f(1)\\ f(0)\ea\right)
=\left(\ba{c}f(0)\\ 0\ea\right),\\[17pt]
\dis J^0f&=&\dis\left(\ba{cc}\ha&0\\ 0&-\ha\ea\right)
\left(\ba{c}f(1)\\ f(0)\ea\right)
=\left(\ba{c}\ha f(1)\\ -\ha f(0)\ea\right).
\ec
We equip $\C^{\{0,1\}}$ and the space in (\ref{tensC}) with the standard
inner product, which has the consequence that $A^\ast=A^\dgg$ and
\[
(J^-_i)^\ast=J^+_i,\quad(J^+_i)^\ast=J^-_i,\quand(J^0_i)^\ast=J^0_i,
\]
showing that the operators $J^\pm_i$ and $J^0_i$ define a unitary
representation of our Lie algebra.

According to the general principle (\ref{prodtwine}), to find an intertwiner
$D$ which acts on the product space (\ref{tensC}), it suffices to find an
intertwiner for the two-dimensional space corresponding to a single site, and
then take the product over all sites. Setting
\[
Q:=\left(\ba{cc}0&1\\ 1&0\ea\right),
\]
it is straightforward to check that
\[
J^\pm Q=QJ^\mp=Q(K^\pm)^\dgg\quand
J^0Q=Q(-J^0)=Q(K^0)^\dgg.
\]
Now, for example, if $S=\{1,2,3\}$ consists of only three sites, then in view
of (\ref{prodtwine})
\[
D:=Q\otimes Q\otimes Q\quad\mbox{satisfies}\quad
J^\pm_iD=D(K^\pm_i)^\dgg\quand J^0_iD=D(K^0_i)^\dgg
\]
$(i=1,2,3)$. In terms of matrix elements, we have $Q(x_i,y_j)=1_{\{x_i\neq
  y_i\}}$ and hence the self-duality function of the symmetric exclusion
process that we have found is
\[
D(x,y)=\prod_{i\in S}1_{\txt\{x_i\neq y_i\}}\qquad\big(x,y\in\{0,1\}^S\big).
\]

\subsection{The symmetric inclusion process}\label{S:SIP}

Let $S$ be a finite set and let $\al:S\to(0,\infty)$ and $q:S\times S\to\half$
be functions such that $q(i,j)=q(j,i)$ and $q(i,i)=0$ for each $i\in S$. By
definition, the \emph{Brownian energy process} or \emph{BEP} with parameters
$\al,q$ is the diffusion process $(Z_t)_{t\geq 0}$ with state space $\half^S$
and generator
\be\label{BEP}
L:=\ha\sum_{i,j\in S}q(i,j)
\big[(\al_jz_i-\al_iz_j)(\dif{z_j}-\dif{z_i})
+z_iz_j(\dif{z_j}-\dif{z_i})^2\big].
\ee
This diffusion has the property that $\sum_iZ_t(i)$ is a preserved
quantity. The drift part of the generator is zero if $z_i=\la\al_i$ for some
$\la>0$. If $z_i/\al_i>z_j/\al_j$, then the drift has the tendency to make
$z_i$ smaller and $z_j$ larger.

In analogy with (\ref{Kidef}), we define operators acting on smooth functions
$f:\half^S\to\R$ by:
\be\label{Kidefi}\bac
\dis\Ki^-_if(z)&=&\dis z_i\diff{z_i}f(z)+\al_i\dif{z_i}f(z),\\[5pt]
\dis\Ki^+_if(z)&=&\dis z_if(z),\\[5pt]
\dis\Ki^0_if(z)&=&\dis z_i\dif{z_i}f(z)+\ha\al_if(z).
\ec
By (\ref{Kcomrel}), these operators satisfy the commutation relations
\[
[\Ki^0_i,\Ki^\pm_j]=\pm\de_{ij}\Ki^\pm_i
\quand[\Ki^-_i,\Ki^+_j]=2\de_{ij}\Ki^0_i.
\]
It follows that these operators define a representation of the Lie algebra
\[
\bigoplus_{i\in S}\gk_i,
\]
where each $\gk_i$ is a copy of $\suk(1,1)$, on the product space
\[
\C^{\half^S}\cong(\C^\half)^{\otimes S},
\]
which is the tensor product of $|S|$ copies of $\C^\half$.

We can express the generator (\ref{BEP}) of the Brownian energy process in
terms of the operators from (\ref{Kidefi}) as
\be\label{Lexpres}
L=\ha\sum_{i,j\in S}q(i,j)
\big[\Ki^+_i\Ki^-_j+\Ki^-_i\Ki^+_j-2\Ki^0_i\Ki^0_j+\ha\al_i\al_j\big].
\ee
Note that this is very similar to the expression for the symmetric exclusion
process in (\ref{excl1}).

We define operators acting on functions $f:\N^S\to\R$ by
\be\bac\label{Kdefi}
\dis K^-_if(x)&=&\dis x_if(x-\de_i),\\[5pt]
\dis K^+_if(x)&=&\dis(\al_i+x_i)f(x+\de_i),\\[5pt]
\dis K^0_if(x)&=&\dis(\ha\al_i+x_i)f(x).
\ec
In view of (\ref{ducomrel}), these operators define a representation of the
conjugate of our Lie algebra. It turns out that the conjugate of this
representation is equivalent to the representation defined by the operators in
(\ref{Kidefi}). This is a nontrivial statement that depends crucially on the
fact that the parameters $\al_i$ are the same in both expressions.
Indeed, we have seen in Subsection~\ref{S:SU11} that $\al$ is twice the
  Bargmann index and that representations with a different Bargmann index have
  a different Casimir operator and hence are not equivalent. Letting $\Phi$
denote the intertwiner of $\Ki^\pm_i$ and $(K^\pm_i)^\dgg$, we can write
$\Phi$ in the form (\ref{Phikern}), where by Proposition~\ref{P:Hilb} $D$ is a
duality function. Similar to what we did at the end of Subsection~\ref{S:SEP},
we will choose a duality function of product form:
\be\label{Dprod}
D(z,x)=\prod_{i\in S}Q(z_i,x_i)\qquad(z\in\half^S,\ x\in\N^S),
\ee
where $Q$ is a duality function for the single-site operators, i.e.,
\be
\Ki^\pm Q(\,\cdot\,,x)(z)=K^\pm Q(z,\,\cdot\,)(x),\quad
\Ki^0Q(\,\cdot\,,x)(z)=K^0Q(z,\,\cdot\,)(x)
\ee
$(z\in\half,\ x\in\N)$. It turns out that
\be\label{DGa}
Q(z,x):=\frac{\Ga(\al+x)}{\Ga(\al)}z^x=z^x\prod_{k=0}^{x-1}(\al+k).
\ee
does the trick. This may look a bit complicated but the form of this duality
function can in fact quite easily be guessed from the inductive relation
\[
zQ(z,x)=\Ki^+Q(\,\cdot\,,x)(z)=K^+Q(z,\,\cdot\,)(x)=(\al+x)Q(z,x+1).
\]

Our calculations so far imply that the generator in (\ref{Lexpres}) is dual
with respect to the duality function in (\ref{Dprod})--(\ref{DGa}) to the
generator
\be\label{hatLexpres}
\hat L=\ha\sum_{i,j\in S}q(i,j)
\big[K^-_jK^+_i+K^+_jK^-_i-2K^0_jK^0_i+\ha\al_j\al_i\big].
\ee
It turns out that we are lucky in the sense that this is a Markov generator.
In view of the similarity with (\ref{excl1}) (with the role of $\suk(2)$
replaced by $\suk(1,1)$), in \cite{GRV10}, the corresponding process has been
called the \emph{symmetric inclusion process} or \emph{SIP}. The fact that
$\hat L$ is a Markov generator can be seen by rewriting it as
\be\ba{r@{\,}l}\label{SIP}
\dis\hat L:=\sum_{i,j\in S}q(i,j)\Big[
&\dis\al_jx_i
\big\{f\big(x-\de_i+\de_j\big)-f\big(x\big)\big\}\\[5pt]
&\dis+x_ix_j
\big\{f\big(x-\de_i+\de_j\big)-f\big(x\big)\big\}\Big].
\ec
The Markov process $(X_t)_{t\geq 0}$ with generator $\hat L$ has the property
that $\sum_iX_t(i)$ is a preserved quantity. The terms in the generator
involving the constants $\al_j$ describe a system of independent random walks,
where each particle at $i$ jumps with rate $\al_j$ to the site $j$. A
reversible law for this part of the dynamics is a Poisson field with local
intensity $\la\al_i$ for some $\la>0$. The remaining terms in the generator
describe a dynamics where particles at $i$ jump to $j$ with a rate that is
proportional to the number $x(j)$ of particles at $j$. This part of the
dynamics causes an attraction between particles.

\subsection{Duality functions of product form}\label{S:sudbur}

In the previous two subsections, we have seen that for a Markov process whose
state space is a Carthesian product of other spaces, it is often natural to
choose duality functions of \emph{product form} as in (\ref{Dprod}). This idea
does not depend on Lie algebras and is in fact older than the use of Lie
algebras in duality theory.

In a series of papers \cite{LS95,LS97,Sud00}, Lloyd and Sudbury have
systematically searched for dualities in a large class of interacting particle
systems, which contains many well-known systems such as the voter model,
contact process, and symmetric exclusion process. Let $S$ be a finite set and
let $q:S^2\to\half$ be a function such that $q(i,j)=q(j,i)$ and $q(i,i)=0$ for
all $i\in S$. Let $L=L(a,b,c,d,e)$ be the Markov generator, acting on functions
$f:\{0,1\}^S\to\R$, as
\be\ba{r@{\,}l}\label{Labcde}
\dis Lf(x)=\sum_{i,j\in S}q(i,j)\Big[
&\dis\ha a1_{\{(x(i),x(j))=(1,1)\}}\big\{f(x-\de_i-\de_j)-f(x)\big\}\\[5pt]
&\dis b1_{\{(x(i),x(j))=(0,1)\}}\big\{f(x+\de_i)-f(x)\big\}\\[5pt]
&\dis c1_{\{(x(i),x(j))=(1,1)\}}\big\{f(x-\de_i)-f(x)\big\}\\[5pt]
&\dis d1_{\{(x(i),x(j))=(0,1)\}}\big\{f(x-\de_j)-f(x)\big\}\\[5pt]
&\dis e1_{\{(x(i),x(j))=(0,1)\}}\big\{f(x+\de_i-\de_j)-f(x)\big\}\Big].
\ec
The dynamics of the Markov process with generator $L$ can be described by
saying that for each pair of sites $i,j$, the configuration of the process at
these sites makes the following transitions with the following rates:
\[\ba{r@{\quad}c@{\quad}l@{\quad}l}
11\mapsto 00&\mbox{with rate}& aq(i,j)&\mbox{(annihilation)},\\[5pt]
01\mapsto 11&\mbox{with rate}& bq(i,j)&\mbox{(branching)},\\[5pt]
11\mapsto 01&\mbox{with rate}& cq(i,j)&\mbox{(coalecence)},\\[5pt]
01\mapsto 00&\mbox{with rate}& dq(i,j)&\mbox{(death)},\\[5pt]
01\mapsto 10&\mbox{with rate}& eq(i,j)&\mbox{(exclusion dynamics)}.
\ea\]
Note that the factor $\ha$ in front of $a$ disappears since the total rate of
this transition is $\ha a(q(i,j)+q(j,i))=aq(i,j)$. A lot of well-known
interacting particle systems fall into this class. For example
\[\ba{r@{\quad}l}
\mbox{voter model}&b=d=1,\mbox{ other parameters }0,\\[5pt]
\mbox{contact process}&b=\la,\ c=d=1,\mbox{ other parameters }0,\\[5pt]
\mbox{symmetric exclusion}&e=1,\mbox{ other parameters }0.
\ea\]

As we have already seen in (\ref{tensC}), the class of all
functions $f:\{0,1\}^S\to\R$ can be written as the tensor product
\[
\R^{\,\txt\{0,1\}^S}\cong\bigotimes_{i\in S}\R^{\,\txt\{0,1\}},
\]
with one `factor' $\R^{\,\txt\{0,1\}}$ for each site $i\in S$. Moreover,
duality functions $D$ on the space $\{0,1\}^S\times\{0,1\}^S$ can be viewed as
matrices corresponding to linear operators that act on $\R^{\{0,1\}^S}$.
Based on various arguments that are not very important at this point, Lloyd
and Sudbury decided to look for duality functions of \emph{product form}
\be\label{profo}
D(x,y)=\prod_{i\in S}Q(x_i,y_i),
\ee
where $Q$ is a $2\times 2$ matrix.
After a more or less systematic search for suitable matrices $Q$, Lloyd and
Sudbury find a rich class of dualities for matrices of the form
\be\label{qpsi}
\left(\ba{cc}
Q_q(0,0)&Q_q(0,1)\\
Q_q(1,0)&Q_q(1,1)
\ea\right)
=\left(\ba{cc}
1&1\\
1&q\ea\right),
\ee
where $q\in\R\beh\{1\}$ is a constant. This choice of $Q$ yields the duality
function
\be\label{Dq}
D_q(x,y):=\prod_{i\in S}Q_q(x_i,y_i)
=q^{\,\txt\sum_{i\in S}x_iy_j}\qquad\big(x,y\in\{0,1\}^S\big).
\ee
In particular, setting $q=0$ yields
\[
D_0(x,y)=1_{\{\sum_{i\in S}x_iy_j=0\}},
\]
which corresponds to the well-known \emph{additive systems duality}, while
$q=-1$ is known as \emph{cancellative systems duality}. For these special
values of $q$, the duality can in fact be upgraded to a pathwise duality as in
Subsection~\ref{S:path}, using a construction in terms of open paths in a
graphical representation. Interestingly, for other values of $q$, there seems
to be no pathwise interpretation of the duality with duality function $D_q$.

We cite the following theorem from \cite{LS95,Sud00}. A somewhat more general
version of this theorem which drops the symmetry assumption
$q(i,j)=q(j,i)$ at the cost of replacing (\ref{dupar}) by a somewhat more
complicated set of conditions can be found in \cite[Appendix~A in the version
  on the ArXiv]{Swa06}.

\bt[q-duality]\hspace{-4pt}
The\label{T:qdual} generators $L(a,b,c,d,e)$ and $L(a',b',c',d',e')$ from
(\ref{Labcde}) are dual with respect to the duality function $D_q$ from
(\ref{Dq}) if and only if
\be\label{dupar}
a'=a+2q\ga,\quad b'=b+\ga,\quad c'=c-(1+q)\ga,\quad d'=d+\ga,
\quad e'=e-\ga,
\ee
where $\ga:=(a+c-d+q b)/(1-q)$.
\et

\subsection{Intertwining and thinning}\label{S:twine}

In Subsection~\ref{S:dual}, when we introduced Markov process duality, we also
defined the similar concept of intertwining of Markov processes. So far, we
have not discussed this second concept very much, but it turns out that the two
are closely related. In particular, as Lloyd and Sudbury already observed
\cite{LS95,Sud00}, there is a close connection between q-duality and thinning
relations. To explain this, we start with a general principle, that says that
if two Markov processes are both dual to a third Markov process, then we can
expect an intertwining relation between the first two processes.

\bl[Duality and intertwining]
Let\label{L:doubdu} $\Om$ and $\hat\Om$ be finite sets, and let
$L_i:\R^\Om\to\R^\Om$, $\hat L:\R^{\hat\Om}\to\R^{\hat\Om}$, and
$D_i:\R^{\hat\Om}\to\R^\Om$ be linear operators such that
\be
L_iD_i=D_i\hat L^\dgg\qquad(i=1,2).
\ee
Assume that $D_1$ and $D_2$ are invertible. Then
\be\label{twodu}
L_1(D_1D_2^{-1})=(D_1D_2^{-1})L_2.
\ee
\el
\bpro
This follows by writing $D_1^{-1}L_1D_1=\hat L^\dgg=D_2^{-1}L_2D_2$.
\epro

We have seen that for interacting particle systems, there are good reasons to
look for duality functions of product form as in (\ref{profo}). Likewise, it
is natural to look for intertwining probability kernels of product form. If
the state space is of the form $\{0,1\}^S$, this means that we are looking for
kernels of the form
\[
K(x,y)=\prod_{i\in S}M(x_i,y_i)\qquad\big(x,y\in\{0,1\}^S\big),
\]
where $M$ is a probability kernel on $\{0,1\}$. If we moreover require that
$M(0,0)=1$ (which is natural for interacting particle systems for which the
all zero state is a trap), then there is only a one-parameter family of such
kernels. For $p\in[0,1]$, let $M_p$ be the probability kernel on $\{0,1\}$
given by
\be\label{Mp}
M_p=\left(\ba{cc}M_p(0,0)&M_p(0,1)\\
M_p(1,0)&M_p(1,1)\ea\right)
:=\left(\ba{cc}1&0\\
1-p&p\ea\right),
\ee
and let
\be\label{Kp}
K_p(x,y):=\prod_{i\in S}M_p(x_i,y_i)\qquad\big(x,y\in\{0,1\}^S\big)
\ee
the corresponding kernel on $\{0,1\}^S$ of product form. We can interpret a
configuration of particles, where $x_i=1$ if the site $i$ is occupied by a
particle, and $x_i=0$ otherwise. Then $K_p$ is a \emph{thinning} kernel that
independently for each site throws away particles with probability $1-p$ or
keeps them with probability $p$. It is easy to see that
\[
K_pK_{p'}=K_{pp'},
\]
i.e., first thinning with $p$ and then with $p'$ is the same as thinning with
$pp'$. There is a close relation between Lloyd and Sudbury's duality function
$D_q$ from (\ref{Dq}) and thinning kernels of the form (\ref{Kp}). We claim
that
\be\label{DqDq}
D_qD_{q'}^{-1}=K_p\quad\mbox{with}\quad p=\frac{1-q}{1-q'}
\qquad(q,q'\in\R,\ q'\neq 1).
\ee
Since both $D_q$ and $K_p$ are of product form, i.e.,
\[
D_q=\bigotimes_{i\in S}Q_q
\quand
K_p=\bigotimes_{i\in S}M_p
\]
with $Q_q$ and $M_p$ as in (\ref{qpsi}) and (\ref{Mp}), it suffices to check
that
\[
Q_qQ_{q'}^{-1}=M_p\quad\mbox{with}\quad p=\frac{1-q}{1-q'}.
\]
Indeed, one can check that
\[
Q_q^{-1}=\left(\ba{cc}
1&1\\
1&q
\ea\right)^{-1}
=(1-q)^{-1}\left(\ba{cc}
-q&1\\
1&-1
\ea\right)\qquad(q\neq 1),
\]
and that
\[
Q_qQ_{q'}^{-1}=(1-q')^{-1}
\left(\ba{cc}
1&1\\
1&q
\ea\right)
\left(\ba{cc}
-q'&1\\
1&-1
\ea\right)
=\left(\ba{cc}
1&0\\
\frac{q-q'}{1-q'}&\frac{1-q}{1-q'}
\ea\right)=M_p,
\]
as claimed.

\bp[Thinning and $q$-duality]
Let\label{P:thinq} $L_1$ and $L_2$ be generators of Markov processes with
state space $\{0,1\}^S$. Assume that there exists an operator $\hat L$ such that
\be
L_iD_{q_i}=D_{q_i}\hat L^\dgg\qquad(i=1,2)
\ee
for some $q_1,q_2\in\R$ such that $q_2\neq 1$ and
$p:=(1-q_1)/(1-q_2)\in[0,1]$. Then
\be
L_1K_p=K_pL_2.
\ee
\ep

\bpro
This follows from (\ref{DqDq}) and Lemma~\ref{L:doubdu}. Note that in general,
there is no guarantee that the operator $D_1D_2^{-1}$ from
Lemma~\ref{L:doubdu} is a probability kernel. In a way,
Proposition~\ref{P:thinq} explains why the q-duality function $D_q$ is
natural, because it is closely linked to the natural concept of thinning.
\epro

\subsection{The biased voter model}\label{S:bias}

In this section, we demonstrate Lloyd-Sudbury theory on the example of the
\emph{biased voter model} with \emph{selection parameter} $s>0$, which is the
interacting particle system with generator
\[
L(a,b,c,d,e)=L(0,1+s,0,1,0)=:L_{\rm bias}.
\]
We apply Theorem~\ref{T:qdual} to find $q$-duals of the biased voter model.
For simplicity, we restrict ourselves here to dual generators of the form
$L(a',b',c',d',e')$ with $a'=0$, which means that we must choose the parameter
$q$ as
\[
q=0\quad\mbox{or}\quad q=(1+s)^{-1}.
\]
For $q=0$ we find the dual generator
\[
L(a',b',c',d',e')=L(0,s,1,0,1)=:L_{\rm braco},
\]
which describes a system of branching and coalescing random walks with
branching parameter $s$. For $q=(1+s)^{-1}$, we find a self-duality, i.e.,
in this case $L(a',b',c',d',e')=L(a,b,c,d,e)=L_{\rm bias}$.

Since $L_{\rm bias}$ and $L_{\rm braco}$ are both $q$-dual to $\hat L=L_{\rm
  bias}$, Proposition~\ref{P:thinq} tells us that there is a thinning relation
between biased voter models and systems of branching and coalescing random
walks of the form
\[
L_{\rm bias}K_p=K_pL_{\rm braco}
\quad\mbox{with}\quad
p=\frac{1-(1+s)^{-1}}{1-0}=\frac{s}{1+s}.
\]
As explained in Subsection~\ref{S:twine}, this implies that if we start a
biased voter model $(X_t)_{t\geq 0}$ and a system of branching and coalescing
random walks $(Y_t)_{t\geq 0}$ in initial states $\mu^{\rm bias}_t$ and
$\mu^{\rm braco}_t$ denote the laws of $X_t$ and $Y_t$, then
\[
\mu^{\rm bias}_0K_p=\mu^{\rm braco}_0
\quad\mbox{implies}\quad
\mu^{\rm bias}_tK_p=\mu^{\rm braco}_t\quad(t\geq 0).
\]
In other words, the following two procedures are equivalent:
\begin{enumerate}
\item Evolve a particle configuration for time $t$
  according to biased voter model dynamics, then thin with $p$.
\item Thin a particle configuration with $p$, then evolve for time $t$
  according to branching coalescing random walk dynamics.
\end{enumerate}
In particular, if we start $X$ in the initial state $X_0(i)=1$ for all $i\in
S$, then because of the nature of the voter model, we will have $X_t(i)=1$ for
all $i\in S$ and $t\geq 0$. Applying the thinning relation now shows that
product measure with intensity $p$ is an invariant law for branching
coalescing random walk dynamics. Thus, there is a close connection between:
\begin{itemize}
\item[I.] $q$-duality,
\item[II.] thinning relations,
\item[III.] invariant laws of product form.
\end{itemize}
Although Lloyd-Sudbury theory is restricted to Markov processes with state
space of the form $\{0,1\}^S$, many other dualities, including the
self-duality of the Wright-Fisher diffusion from Section~\ref{S:WFself}, can
be derived from Lloyd-Sudbury duals by taking a suitable limit \cite{Swa06}.

\subsection{Time-reversal and symmetry}\label{S:symm}

In this subsection we present an idea from \cite{GKRV09}, which says that
nontrivial dualities can sometimes be found by starting from a ``trivial''
duality which is based on time reversal, and then using a symmetry of the
model to transform such a duality into a nontrivial one. Although Lie algebras
are not strictly needed in this approach, writing generators in terms of the
basis elements of a representation of a Lie algebra can help finding suitable
symmetries.

Each irreducible Markov process with finite state space $\Om$ has a unique
invariant measure, i.e., a probability measure $\mu$ such that
\[
\mu L=0\quad\mbox{or equivalently}\quad\mu P_t=\mu\quad(t\geq 0),
\]
where $L$ denotes the generator and $(P_t)_{t\geq 0}$ the semigroup of the
Markov process. Irreducibility implies that $\mu(x)>0$ for all
$x\in\Om$. Letting $(X_t)_{t\in\R}$ denote the stationary process, we see that
the semigroup $(\ti P_t)_{t\geq 0}$ of the time-reversed process is given by
\[\ba{l}
\dis\ti P_t(x,y)=\frac{\P[X_0=y,\ X_t=x]}{\P[X_t=x]}\\[5pt]
\dis\quad=\frac{\mu(y)P_t(y,x)}{\mu(x)}=\mu(y)P_t(y,x)\mu(x)^{-1}\quad(t\geq 0).
\ea\]
Differentiating shows that the generator $\ti L$ of the
time-reversed process is given by\footnote{This formula is wrong in
  \cite[below (12)]{GKRV09}.}
\[
\ti L(x,y)=\mu(y)L(y,x)\mu(x)^{-1}.
\]
Let $R$ denote the diagonal matrix
\[
R(x,y):=\de_{x,y}\mu(x)^{-1}.
\]
Then $L(y,x)\mu(x)^{-1}=\ti L(x,y)\mu(y)^{-1}=\mu(y)^{-1}\ti L^\dgg(y,x)$ can
be rewritten as
\[
LR=R\ti L^\dgg,
\]
which shows that $\ti L$ is dual to $L$ with duality function $R$.
In particular, reversible processes (for which $\ti L=L$) are always
\emph{self-dual} with duality function $R(x,y)$. Note that since $R$ is
diagonal, it is invertible with
\[
R^{-1}(x,y):=\de_{x,y}\mu(x)\qquad(x,y\in\Om).
\]

Let $V$ be a finite dimensional complex linear space and let $L:V\to V$ be any
linear operator (not necessarily a Markov generator). Then it is known that
there exists an invertible matrix $Q\in\Li(V)$ such that
\be\label{LtL}
LQ=QL^\dgg\quad\mbox{or equivalently}\quad
L^\dgg Q^{-1}=Q^{-1}L
\ee
Thus, \emph{every} finite dimensional linear operator is self-dual and the
self-duality function $Q$ can be chosen such that it is invertible, viewed as
a matrix. Let
\[
\Ci_L:=\{A\in\Li(V):AL=LA\}
\]
be the algebra of all elements of $\Li(V)$ that commute with $L$.
We call this the space of \emph{symmetries} of $L$. In
\cite[Thm~2.6]{GKRV09}, the following simple observation is made.

\bl[Self-duality functions]
Let\label{L:selfsym} $L$ be a linear operator on some finite dimensional
linear space $V$. Fix some $Q$ as in (\ref{LtL}). Then the set of all
self-duality functions of $L$ is given by
\[
\{SQ:S\in\Ci_L\}.
\]
\el
\bpro
Clearly, if $S\in\Ci_L$, then
\[
LSQ=SLQ=SQL^\dgg,
\]
showing that $SQ$ is a self-duality function. Conversely, if $D$ is a
self-duality function, then we can write $D=SQ$ with $S=DQ^{-1}$. Now,
since $D$ is a self-duality function,
\[
SL=DQ^{-1}L=DL^\dgg Q^{-1}=LDQ^{-1}=LS,
\]
which shows that $S\in\Ci_L$.
\epro

For dualities, we can play a similar game. Once we have two operators $L,\hat
L$ that are dual with duality function $D$, i.e.,
\[
LD=D\hat L^\dgg,
\]
we have that for any $S\in\Ci_L$, the operators $L,\hat
L$ are also dual with duality function $SD$, as follows by writing
\[
LSD=SLD=SD\hat L^\dgg.
\]
If $D$ is invertible, then every duality function of $L$ and $\hat L$ is of
this form. Indeed, if $\ti D$ is any duality function, then we can write $\ti
D=SD$ with $S=\ti DD^{-1}$. Now
\[
SL=\ti DD^{-1}L=\ti DL^\dgg D^{-1}=L\ti DD^{-1}=LS,
\]
proving that $S\in\Ci_L$. See also \cite[Thm~2.10]{GKRV09}.

\subsection{The symmetric exclusion process revisited}\label{S:SEP2}

Following \cite[Sect.~3.1]{GKRV09}, we demonstrate the principles explained in
the previous subsections to derive a self-duality of the symmetric exclusion
process. Our starting point is formula (\ref{excl1}), which expresses the
generator $L$ in terms of operators $J^\pm_i,J^0_i$ that define a
representation $(V,\pi)$ of a Lie algebra $\gk$ that is the direct sum of
finitely many copies of the Lie algebra $\suk(2)$, with one copy for each site
$i\in S$. Since $r(i,j)=r(j,i)$, we can rewrite this formula as
\be
L=\ha\sum_{i,j}r(i,j)\big[J^-_iJ^+_j+J^-_jJ^+_i+2J^0_iJ^0_j-\ha I\big].
\ee
A straightforward calculation shows that
\be\label{Lsym}
\sum_k[J^\pm_k,L]=0\quand\sum_k[J^0_k,L]=0\qquad(k\in S).
\ee

We need a bit of general theory. If $U,V,W$ are representations of the same
Lie algebra $\gk$, then we can equip their tensor product $U\otimes V\otimes
W$ with the structure of a representation of $\gk$ by putting
\be\label{tensrep}
A(u\otimes v\otimes w):=
Au\otimes v\otimes w+u\otimes Av\otimes w+u\otimes v\otimes Aw
\qquad(A\in\gk),
\ee
and similar for the tensor product of any finite number of representations,
see formula (\ref{altens}) in the appendix. This definition also naturally
equips $U\otimes V\otimes W$ with the structure of a representation of the Lie
group $G$ associated with $\gk$, in such a way that
\[
e^{tA}(u\otimes v\otimes w)=e^{tA}u\otimes e^{tA}v\otimes e^{tA}w
\qquad(A\in\gk,\ t\geq 0),
\]
where for each $A\in\gk$, the operator $e^{tA}$ is an element of the Lie group
$G$ associated with $\gk$. Thus, the representation (\ref{tensrep})
corresponds to letting the Lie group act in the same way on each space in the
tensor product.

In our specific set-up, this means that the operators $K^-,K^+,K^0$ defined by
\be
K^-:=\sum_kJ^-_k,\quad K^+:=\sum_kJ^+_k,\quad K^0:=\sum_kJ^0_k
\ee
define a representation of $\suk(2)$ on the product space
\[
\C^{\{0,1\}^S}\cong\bigotimes_{i\in S}\C^{\{0,1\}}.
\]
(Indeed, one can check that $K^-,K^+,K^0$ satisfy the commutation relations of
$\suk(2)$.) Let $c_-K^-+c_+K^++c_0K^0$ be an operator in the linear space
spanned by $K^-,K^+,K^0$. Then
\be\label{prodsym}
\ex{t(c_-K^-+c_+K^++c_0K^0)}=\bigotimes_{i\in S}\ex{t(c_-J^-+c_+J^++c_0J^0)}
\qquad(t\geq 0),
\ee
i.e., a natural group of symmetries of the generator $L$ is formed by all
operators of the form (\ref{prodsym}) and their products, and this actually
corresponds to a representation of the Lie group ${\rm SU}(2)$.

We take this as our motivation to look at one specific operator of the form 
(\ref{prodsym}), which is $e^{K^+}$. One can check that the uniform
distribution is an invariant law for the exclusion process, so by the
principle of Subsection~\ref{S:symm}, the function
\[
D(x,y)=1_{\{x=y\}}=\prod_{i\in S}1_{\{x_i=y_i\}}
\]
is a trivial self-duality function. Applying Lemma~\ref{L:selfsym} to the
symmetry $S=e^{K^+}$, we see that $SD=SI=S$ is also a self-duality function.
Since $S$ factorizes over the sites, it suffices to calculate $S$ for a single
site, and then take the product. We recall from (\ref{Jpmo}) that
\[
J^+f\left(\ba{cc}0&1\\ 0&0\ea\right)
\left(\ba{c}f(1)\\ f(0)\ea\right)
=\left(\ba{c}f(0)\\ 0\ea\right),
\]
which gives
\[
\ex{J^+}=\sum_{n=0}^\infty\frac{1}{n!}(J^+)^n
=I+J^+=\left(\ba{cc}1&1\\ 0&1\ea\right)
\]
and finally yields the duality function
\[
S(x,y)=\prod_{i\in S}1_{\txt\{x_i\geq y_i\}}\qquad\big(x,y\in\{0,1\}^S\big).
\]

\appendix

\section{A crash course in Lie algebras}\label{A:Lie}

\subsection{Lie groups}

In the present appendix, we give a bit more background on Lie algebras. In
particular, we explain how Lie algebras are closely linked to Lie groups, and
how every Lie algebra can naturally be embedded in an algebra, called the
universal enveloping algebra. We also explain how properties of the Lie group
(in particular, compactness) are related to representations of its associated
Lie algebra.

A \emph{group} is a set $G$ which contains a special element $I$, called the
\emph{identity}, and on which a \emph{group product} $(A,B)\mapsto AB$ and
\emph{inverse operation} $A\mapsto A^{-1}$ are defined such that
\begin{enumerate}
\item $IA=AI=A$
\item $(AB)C=A(BC)$
\item $A^{-1}A=AA^{-1}=I$.
\end{enumerate}
A group is \emph{abelian} (also called \emph{commutative}) if $AB=BA$ for all
$A,B\in G$. A \emph{group homomorphism} is a map $\Phi$ from one group $G$
into another group $H$ that preserves the group structure, i.e.,
\begin{enumerate}
\item $\Phi(I)=I$,
\item $\Phi(AB)=\Phi(A)\Phi(B)$,
\item $\Phi(A^{-1})=\Phi(A)^{-1}$.
\end{enumerate}
If $\Phi$ is a bijection, then $\Phi^{-1}$ is also a group homomorphism. In this
case, we call $\Phi$ a \emph{group isomorphism}. A \emph{subgroup} of a group
$G$ is a subset $H\sub G$ such that $I\in H$ and $H$ is closed under the
product and inverse, i.e., $A,B\in H$ imply $AB\in H$ and $A\in H$ implies
$A^{-1}\in H$. A subgroup is in a natural way itself a group.

A \emph{Lie group} is a smooth manifold $G$ which is also a group such that
the group product and inverse functions
\[
G\times G\ni(A,B)\mapsto AB\in G
\quand
G\ni A\mapsto A^{-1}\in G
\]
are smooth. A \emph{finite-dimensional representation} of $G$ is a
finite-dimensional linear space $V$ over $\R$ or $\C$ together with a map
\[
G\times V\ni(A,v)\mapsto Av\in V
\]
such that
\begin{enumerate}
\item $v\mapsto Av$ is linear,
\item $Iv=v$,
\item $A(Bv)=(AB)v$.
\end{enumerate}
Letting $\Li(V)$ denote the space of all linear operators $A:V\to V$,
these conditions are equivalent to saying that the map $\Pi:G\to\Li(V)$
defined by
\[
\Pi(A)v:=Av
\]
is a group homomorphism from $G$ into the \emph{general linear group} ${\rm
  GL}(V)$ of all invertible linear maps $A:V\to V$. A representation is
\emph{faithful} if $\Pi$ is one-to-one, i.e., if $A\mapsto\Pi(A)$ is a group
isomorphism between $G$ and the subgroup $\Pi(G):=\{\Pi(A):A\in G\}$ of ${\rm
  GL}(V)$.

One can prove that if $G$ is a Lie group and $V$ is a faithful
finite-dimensional representation, then $\Pi(G)$ is a closed subset of ${\rm
  GL}(V)$ and $\Pi:G\to\Pi(G)$ is a homeomorhism. Conversely, each
closed subgroup of ${\rm GL}(V)$ is a Lie group.  Such Lie groups are called
\emph{matrix Lie groups}. Not every Lie group has a finite dimensional
faithful representation, so not every Lie group is a matrix Lie group, but
many important Lie groups are matrix Lie groups and following \cite{Hal03} we
will mostly focus on them from now on.

\subsection{Lie algebras}

An \emph{algebra} is a finite-dimensional linear space $\ak$ over $\R$ or
$\C$ with a special element $I$ called \emph{unit element}
and on which there is defined a product
\[
\ak\times\ak\ni(A,B)\mapsto AB\in\ak
\]
such that
\begin{enumerate}
\item $(A,B)\mapsto AB$ is bilinear,
\item $IA=AI=A$,
\item $(AB)C=A(BC)$.
\end{enumerate}
In some textbooks, algebras are not required to contain a unit element.  We
speak of a \emph{real} resp.\ \emph{complex} algebra depending on whether
$\ak$ is a linear space over $\R$ or $\C$. An algebra is \emph{abelian} if
$AB=BA$ for all $A,B\in G$. In any algebra, the
\emph{commutator} of two elements $A,B$ is defined as $[A,B]=AB-BA$. 
If $V$ is a linear space, then $\Li(V)$ is an algebra.

An \emph{algebra homomorphism} is a map $\phi:\ak\to\bk$ from one algebra into
another that preserves the structure, i.e.,
\begin{enumerate}
\item $\phi$ is linear,
\item $\phi(I)=I$,
\item $\phi(AB)=\phi(A)\phi(B)$.
\end{enumerate}
Algebra homomorphisms that are bijections have the property that $\phi^{-1}$
is also a homomorphism; these are called algebra isomorphisms. A
\emph{subalgebra} of an algebra $\ak$ is a linear subspace $\bk\sub\ak$ that
contains $I$ and is closed under the product.

\emph{Lie algebras}, \emph{Lie algebra homomorphisms},
and \emph{isomorphisms} have already been defined in
Section~\ref{S:Lieint}. A \emph{sub-Lie-algebra} is a linear subspace
$\hk\sub\gk$ such that
\[
A,B\in\hk\quad\mbox{implies}\quad[A,B]\in\hk.
\]
If $\gk$ is an algebra, then $\gk$, equipped with the commutator map
$[\,\cdot\,,\,\cdot\,]$, is a Lie algebra. As the example in
Section~\ref{S:Lieint} shows. Lie algebras need not be an algebras.

A \emph{representation} of an algebra $\ak$ is a linear space $V$ together
with a map $\ak\times V\to V$ that satisfies
\begin{enumerate}
\item $(A,v)\mapsto Av$ is bilinear,
\item $Iv=v$,
\item $A(Bv)=(AB)v$.
\end{enumerate}
If $\ak$ is a complex algebra, then we require $V$ to be a linear space over
$\C$, but even when $\ak$ is a real algebra, it is often useful to allow for
the case that $V$ is a linear space over $\C$. In this case, bilinearity means
real linearity in the first argument and complex linearity in the second
argument. We speak of \emph{real} or \emph{complex} representations depending
on whether $V$ is a linear space over $\R$ or $\C$.

A representation $V$ of an algebra $\ak$ gives in a natural way rise to an
algebra homomorphism $\pi:\ak\to\Li(V)$ defined as
\[
\pi(A)v:=Av\qquad(A\in\ak,\ v\in V).
\]
Conversely, given an algebra homomorphism $\pi:\ak\to\Li(V)$ we can equip $V$
with the structure of a representation by defining $Av:=\pi(A)v$. Thus, a
representation $V$ of an algebra $\ak$ is equivalent to a pair $(V,\pi)$ where
$V$ is a linear space and $\pi:\ak\to\Li(V)$ is an algebra homomorphism.  A
representation $(V,\pi)$ is \emph{faithful} if $\pi$ is an isomorphism between
$\ak$ and the subalgebra $\pi(\ak)=\{\pi(A):A\in\ak\}$ of $\Li(V)$.

Representations of Lie algebras have already been defined in
Section~\ref{S:rep}. If $V$ is a complex representation of a real algebra or
Lie algebra $\ak$, then the image of $\ak$ under $\pi$ is only a real subspace
of $\Li(V)$. We can define a complex algebra or Lie algebra $\ak_\C$ whose
elements can formally be written as $A+iB$ with $A,B\in\ak$; this is called
the \emph{complexification} of $\ak$. Then $\pi$ extends uniquely to a
homomorphism from $\ak_\C$ to $\Li(V)$, see \cite[Prop.~3.39]{Hal03}, so $V$
is also a representation of $\ak_\C$.

Every algebra has a faithful representation. Indeed, $\ak$ together with the
map $(A,B)\mapsto AB$ is a representation of itself, and it is not hard to see
(using our assumption that $I\in\ak$) that this representation is faithful.
Lie algebras can be represented on themselves in a construction that
is very similar to the one for algebras.

\bl[Lie algebra represented on itself]
A\label{L:Lierep} Lie algebra $\gk$, equipped with the map
$(A,B)\mapsto[A,B]$, is a representation of itself.
\el
\bpro
It will be convenient to use somewhat different notation for the Lie bracket.
If $\gk$ is a Lie algebra and $X\in\gk$, then we define ${\rm ad}_X:\gk\to\gk$
by
\[
{\rm ad}_X(A):=[X,A].
\]
We need to show that $\gk\ni X\mapsto{\rm ad}_X\in\Li(\gk)$ is a Lie algebra
homomorphism. Bilinearity follows immediately from the bilinear property (i)
of the Lie bracket, so it remains to show that
\[
{\rm ad}_{[X,Y]}(Z)={\rm ad}_X({\rm ad}_Y(Z))-{\rm ad}_Y({\rm ad}_X(Z)).
\]
This can be rewritten as
\[
[[X,Y],Z]=[X,[Y,Z]]-[Y,[X,Z]].
\]
Using also the skew symmetric property~(ii) of the Lie bracket, this can be
rewritten as
\[
0=[Z,[X,Y]]+[X,[Y,Z]]+[Y,[Z,X]],
\]
which is the Jacobi identity.
\epro

In general, representing a Lie algebra on itself as in Lemma~\ref{L:Lierep}
need not yield a faithful representation. (For example, any abelian algebra is
also a Lie algebra and for such Lie algebras ${\rm ad}_X=0$ for each $X$.)
By definition, the \emph{center} of a Lie algebra $\gk$ is the set
\be\label{center}
\{X\in\gk:[X,A]=0\ \forall A\in\gk\}.
\ee
We say that the center is \emph{trivial} if it contains only the zero
element. If $\gk$ has a trivial center, then the representation $X\mapsto{\rm
  ad}_X$ of $\gk$ on itself is faithful. Indeed, ${\rm ad}_X={\rm ad}_Y$
implies $[X,A]=[Y,A]$ for all $A\in\gk$ and hence $X-Y$ is an element of the
center of $\gk$. If the center is trivial, this implies $X=Y$.

\subsection{Relation between Lie groups and Lie algebras}\label{S:gralg}

Let $V$ be a linear space and let $G\sub{\rm GL}(V)$ be a matrix Lie group. 
By definition, the \emph{Lie algebra} $\gk$ of $G$ is the space of all
matrices $A$ such that there exists a smooth curve $\ga$ in $G$ with
\[
\ga(0)=I\quand\dif{t}\ga(t)\big|_{t=0}=A.
\]
In manifold terminology, this says that $\gk$ is the \emph{tangent space} to
$G$ at $I$. For any matrix $A$, we define
\be\label{expdef}
e^A:=\sum_{k=0}^\infty\frac{1}{n!}A^n.
\ee
The following lemma follows from \cite[Cor.~3.46]{Hal03}. The main idea behind
this lemma is that the elements of the Lie algebra act as ``infinitesimal
generators'' of the Lie group.

\bl[Exponential formula]
Let\label{L:exp} $\gk$ be the Lie algebra of a Lie group $G\sub{\rm
  GL}(V)$. Then the following conditions are equivalent.
\begin{enumerate}
\item $A\in\gk$
\item $e^{tA}\in G$ for all $t\in\R$.
\end{enumerate}
\el

The following lemma (a precise proof of which can be found in \cite[Thm~3.20]{Hal03}) says that our terminology is justified.

\bl[Lie algebra property]
The Lie algebra of any matrix Lie group is a real Lie algebra.
\el
\bpro[(sketch)]
Let $\la\in\R$ and $A\in\gk$. By assumption, there exists a smooth curve $\ga$
such that $\ga(0)=I$ and $\dif{t}\ga(t)\big|_{t=0}=A$. But now
$t\mapsto\ga(\la t)$ is also smooth and $\dif{t}\ga(\la t)\big|_{t=0}=\la A$,
showing that $\gk$ is closed under multiplication with real scalars.

Also, if $A,B\in\gk$, then in the limit as $t\to 0$,
\[
e^{tA}e^{tB}=\big((I+tA+O(t^2)\big)\big((I+tB+O(t^2)\big)
=I+(A+B)t+O(t^2),
\]
which suggests that $A+B$ lies in the tangent space to $G$ at $I$; making this
idea precise proves that indeed $A+B\in\gk$, so $\gk$ is a real linear space.

To complete the proof, we must show that $[A,B]\in\gk$ for all $A,B\in\gk$.
It is easy to see that for any $A,B\in\gk$, as $t\to 0$
\[
[e^{tA},e^{tB}]=t^2[A,B]+O(t^3),
\]
and hence
\[
e^{tA}e^{tB}e^{-tA}e^{-tB}
=e^{tA}\{e^{-tA}e^{tB}+[e^{tB},e^{-tA}]\}e^{-tB}
=I+t^2[A,B]+O(t^3).
\]
Since $e^{tA}e^{tB}e^{-tA}e^{-tB}\in G$, this suggests that $[A,B]$ lies in
the tangent space to $G$ at $I$.
\epro

By \cite[Cor.~3.47]{Hal03}, if $\gk$ is the Lie algebra of a Lie group $G$,
then there exist open neighbourhoods $0\in O\sub\gk$ and $I\in U\sub G$ such
that the map
\[
O\ni A\mapsto\ex{A}\in U
\]
is a homeomorphism (a continuous bijection whose inverse is also continuous).
The \emph{identity component} $G_0$ of a Lie group $G$ is the connected
component that contains the identity. By \cite[Prop.~1.10]{Hal03}, $G_0$ is a
subgroup\footnote{In fact, $G_0$ is a normal subgroup -see formula
  (\ref{normal}) below for the definition of a normal subgroup.} of $G$. If
$U$ is an open neighbourhood of $I$, then each element of $G_0$ can be written
as the product of finitely many elements of $U$. In particular, if $G$ is
connected, then $U$ generates $G$. Therefore (see \cite[Cor.~3.47]{Hal03}), if
$G$ is a connected Lie group, then each element $X\in G$ can be written as
\be\label{connect}
X=e^{A_1}\cdots e^{A_n}
\ee
for some $A_1,\ldots,A_n\in\gk$. As \cite[Example~3.41]{Hal03} shows, even if
$G$ is connected, it is in general not true that for each $A,B\in\gk$ there
exists a $C\in\gk$ such that $e^Ae^B=e^C$ and hence in general
$\{e^A:A\in\gk\}$ need not be a group; in particular, this is not always $G$.

Anyway, the Lie algebra uniquely characterizes the local structure of a Lie
group, so it should be true that if two Lie groups $G$ and $H$ are isomorphic,
then their Lie algebras $\gk$ and $\hk$ are also isomorphic. Indeed, by
\cite[Thm.~3.28]{Hal03}, each Lie group homomorphism $\Phi:G\to H$ gives rise
to a unique homomorphism $\phi:\gk\to\hk$ of Lie algebras such that
\be\label{Phiphi}
\Phi(e^A)=e^{\phi(A)}\qquad(A\in\gk).
\ee
In general, the converse conclusion cannot be drawn, i.e., two different Lie
groups may have the same Lie algebra. By definition, a Lie group $G$ is
\emph{simply connected} if it is connected and ``has no holes'', i.e., every
continuous loop can be continuously shrunk to a point. (E.g., the surface of a
ball is simply connected but a torus is not.) We cite the following theorem
from \cite[Thm.~5.6]{Hal03}.

\bt[Simply connected Lie groups]
Let\label{T:simpcon} $G$ and $H$ be matrix Lie groups with Lie algebras $\gk$
and $\hk$ and let $\phi:\gk\to\hk$ be a homomorphism of Lie algebras. If $G$
is simply connected, then there exists a unique Lie group homomorphism
$\Phi:G\to H$ such that (\ref{Phiphi}) holds.
\et

In particular (\cite[Cor.~5.7]{Hal03}), this implies that two simply connected
Lie groups are isomorphic if and only if their Lie algebras are isomorphic.
Every connected Lie group $G$ has a \emph{universal cover} $(H,\Phi)$
(this is stated without proof in \cite[Sect.~5.8]{Hal03}), which
is a simply connected Lie group $H$ together with a Lie group homomorphism
$\Phi:H\to G$ such that the associated Lie algebra homomorphism as in
(\ref{Phiphi}) is a Lie algebra isomorphism. The following lemma says that
such a universal cover is unique up to natural isomorphisms.

\bl[Uniqueness of the universal cover]
Let\label{L:covuni} $G$ be a connected Lie group and let $(H_i,\Phi_i)$
$(i=1,2)$ be universal covers of $G$. Then there exists a unique Lie group
isomorphism $\Psi:H_1\to H_2$ such that $\Psi(\Phi_1(A))=\Phi_2(A)$ $(A\in G)$.
\el
\bpro
Let $\phi_i:\gk\to\hk_i$ denote the Lie algebra homomorphism associated with
$\Phi_i$ as in (\ref{Phiphi}). If a Lie group isomorphism $\Psi$ as in the
lemma exists, then the associated Lie algebra isomorphism $\psi$ must satisfy
$\psi\circ\phi_1=\phi_2$. By assumption, $\phi_i$ $(i=1,2)$ are
isomorphisms, so setting $\psi:=\phi_2\circ\phi_1^{-1}$ defines a Lie algebra
isomorphism from $\hk_1$ to $\hk_2$. By assumption, $H_1$ is simply
connected, so by Theorem~\ref{T:simpcon}, there exists a unique Lie group
homomorphism $\Psi:H_1\to H_2$ such that $\Psi(e^A)=e^{\psi(A)}$ $(A\in\hk_1)$.
Similarly, there exists a unique Lie group homomorphism $\ti\Psi:H_2\to H_1$
such that $\ti\Psi(e^A)=e^{\psi^{-1}(A)}$ $(A\in\hk_2)$. Now
\[
\ti\Psi(\Psi(e^A))=\ti\Psi(e^{\psi(A)})=e^{\psi^{-1}\circ\psi(A)}=e^A
\qquad(A\in\hk_1)
\]
and similarly $\Psi(\ti\Psi(e^A))$ $(A\in\hk_2)$, which (using the fact that
elements of the form $e^A$ with $A\in\hk_i$ generate $H_i$) proves that $\Psi$
is invertible and $\ti\Psi=\Psi^{-1}$.
\epro

Informally, the universal cover $H$ of $G$ is the unique simply connected Lie
group that has the same Lie algebra as $G$. The universal cover of a matrix
Lie group need in general not be a matrix Lie group. Lie's third theorem
\cite[Thm~5.25]{Hal03} says:

\bt[Lie's third theorem]
Every\label{T:Lie3} real Lie algebra $\gk$ is the Lie algebra of some
connected Lie group $G$.
\et

By \cite[Conclusion~5.26]{Hal03}, we can even take $G$ to be a matrix Lie
group, and by restricting to the identity component we can take $G$ to be
connected. By going to the universal cover, we can also take $G$ to be simply
connected, but in this case we may loose the property that $G$ is a
\emph{matrix} Lie group. Anyway, we can conclude:
\begin{quote}
There is a one-to-one correspondence between Lie algebras and simply connected
Lie groups. Every Lie group has a unique universal cover, which is a simply
connected Lie group with the same Lie algebra.
\end{quote}

Let $G$ be a Lie group with Lie algebra $\gk$ and let $(V,\Pi)$ be a
representation of $G$. Then, by (\ref{Phiphi}), there exists a unique Lie
algebra homomorphism $\pi:\gk\to\Li(V)$ such that
\be\label{Pipi}
\Pi(e^A)=e^{\pi(A)}\qquad(A\in\gk).
\ee
More concretely, one has (see \cite[Prop.~4.4]{Hal03})
\be\label{asrep}
\pi(A)v=\dif{t}\Pi(e^{tA})v\big|_{t=0}\qquad(A\in\gk,\ v\in V).
\ee
We say that $(V,\pi)$ is the representation of $\gk$ \emph{associated} with
the representation $(V,\Pi)$ of $G$. Conversely, if $G$ is simply connected,
then by grace of Theorem~\ref{T:simpcon}, through (\ref{Pipi}), each
representation $(V,\pi)$ of $\gk$ gives rise to a unique associated
representation $(V,\Pi)$ of $G$.

\subsection{Relation between algebras and Lie algebras}\label{A:algLie}

If $\ak$ is an algebra and $\ck\sub\ak$ is any subset of $\ak$, then there
exists a smallest subalgebra $\bk\sub\ak$ such that $\bk$ contains
$\ck$. This algebra consists of the linear span of the unit element $I$ and
all finite products of elements of $\ck$. We call $\bk$ the algebra
\emph{generated} by $\ck$. If $\bk=\ak$, then we say that $\ck$
\emph{generates} $\ak$.

Let $\gk$ be a Lie algebra. By definition, an
\emph{enveloping algebra} for $\gk$ is a pair $(\ak,\ik)$ such that
\begin{enumerate}
\item $\ak$ is an algebra and $\ik:\gk\to\ak$ is a Lie algebra homomorphism.
\item The image $\ik(\gk)$ of $\gk$ under $\ik$ generates $\ak$.
\end{enumerate}
We cite the following theorem from \cite[Thms~9.7 and 9.9]{Hal03}.

\bt[Universal enveloping algebra]
For\label{T:unenv} every Lie algebra $\gk$, there exists an enveloping algebra
$(\ak,\ik)$ with the following properties.
\begin{enumerate}
\item If $(\bk,\jk)$ is an enveloping algebra of $\gk$, then there exists a
  unique algebra homomorphism $\phi:\ak\to\bk$ such that $\phi(\ik(A))=\jk(A)$
  for all $A\in\gk$.
\item If $\{X_1,\ldots,X_n\}$ is a basis for $\gk$, then a basis for $\ak$ is
  formed by all elements of the form
\[
\ik(X_1)^{k_1}\cdots\ik(X_n)^{k_n},
\]
where $k_1,\ldots,k_n\geq 0$ are integers. In particular, these elements are
linearly independent.
\end{enumerate}
\et

An argument similar to the proof of Lemma~\ref{L:covuni} shows that the pair
$(\ak,\ik)$ from Theorem~\ref{T:unenv} is unique up to natural
isomorphisms. We call $(\ak,\ik)$ the \emph{universal enveloping algebra} of
$\gk$ and use the notation $U(\gk):=\ak$. By property~(ii), the map $\ik$ is
one-to-one, so we often identify $\gk$ with its image under $\ik$ and
pretend $\gk$ is a sub-Lie-algebra of $U(\gk)$.

As an immediate consequence of property~(i) of Theorem~\ref{T:unenv}, we
see that if $V$ is a representation of a Lie algebra $\gk$ and
$\pi:\gk\to\Li(V)$ is the associated Lie algebra homomorphism, then
there exists a unique algebra homomorphism $\ov\pi:U(\gk)\to\Li(V)$ such that
$\ov\pi(A)=\pi(A)$ $(A\in\gk)$. (Here we view $\gk$ as a sub-Lie-algebra of
$U(\gk)$.) Conversely, of course, every representation of $U(\gk)$ is also a
representation of $\gk$.

If $(V,\pi)$ is a representation of a Lie algebra $\gk$, then we usually
denote the associated representation of $U(\gk)$ also by $(V,\pi)$, i.e., we
identify the map $\pi$ with its extension $\ov\pi$. Note, however, that a
representation $(V,\pi)$ of a Lie algebra $\gk$ can be faithful even when the
associated representation $(V,\pi)$ of $U(\gk)$ is not. Indeed, by
property~(ii) of Theorem~\ref{T:unenv}, $U(\gk)$ is always infinite
dimensional, even though $\gk$ is finite dimensional, so finite-dimensional
faithful representations of $\gk$ are not faithful when viewed as a
representation of $U(\gk)$.

\subsection{Adjoints and unitary representations}

Let $V$ be a finite dimensional linear space equipped with an inner product
$\li\,\cdot\,|\,\cdot\,\re$, which for linear spaces over $\C$
is conjugate linear in its first argument and linear in its second argument.
Each $A\in\Li(V)$ has a unique \emph{adjoint} $A^\ast\in\Li(V)$ such that
\be\label{adj}
\li A^\ast v|w\re=\li v|Aw\re\qquad(v,w\in V).
\ee
An operator $A$ is \emph{self-adjoint} (also called \emph{hermitian}) if
$A^\ast=A$ and \emph{skew symmetric} if $A^\ast=-A$. A \emph{positive
  operator} is an operator such that $\li v|Av\re\geq 0$ for all $v$.
If $V,W$ are linear spaces equipped with inner products, then an operator
$U\in\Li(V,W)$ is called \emph{unitary} if it preserves the inner product, i.e.,
\be\label{Udef}
\li Uv|Uw\re=\li v|w\re\qquad(v,w\in V).
\ee
In particular, an operator $U\in\Li(V)$ is unitary if and only if it is
invertible and $U^{-1}=U$. If $V$ is a finite dimensional linear space over
$\C$, then for $v\in V$ we define operators $\li v|\in\Li(V,\C)$ and
$|v\re\in\Li(\C,V)$ by
\[
\li v|w:=\li v|w\re\quand|v\re c:=cv.
\]
Then $\li v||w\re$ is an operator in $\Li(\C,\C)$ which we can identify with
the complex number $\li v|w\re$. Moreover, $|v\re\li w|$ is an operator in
$\Li(V)$. An \emph{orthonormal} basis $\{e(1),\ldots,e(n)\}$ of $V$ is a basis
such that $\li e(i)|e(j)\re=\de_{ij}$. Then
\[
A=\sum_{ij}A_{ij}|e(i)\re\li e(j)|,
\]
where $A_{ij}$ denotes the matrix of $A$ with respect to the orthonormal basis 
$\{e(1),\ldots,e(n)\}$. An operator $A\in\Li(V)$ is \emph{normal} if
$[A,A^\ast]=0$. An operator is normal if and only if it is diagonal w.r.t.\ some
orthonormal basis, i.e., if it can be written as
\[
A=\sum_i\la_i|e(i)\re\li e(i)|,
\]
where the $\la_i$ are the eigenvalues of $A$. For operators, the following
properties are equivalent.

\begin{tabular}{r@{\quad}c@{\quad}l}
$A$ is hermitian &$\desd$& $A$ is normal with real eigenvalues,\\
$A$ is skew symmetric &$\desd$& $A$ is normal with
 imaginary eigenvalues,\\
$A$ is positive &$\desd$& $A$ is normal with nonnegative eigenvalues,\\
$A$ is unitary &$\desd$& $A$ is normal with eigenvalues of norm 1.\\[3pt]
\end{tabular}

By definition, a \emph{unitary} representation of a Lie group $G$ is a complex
representation $(V,\Pi)$ where $V$ is equipped with an inner product such that
$\Pi(A)$ is a unitary operator for each $A\in G$. A \emph{unitary}
representation of a real Lie algebra $\gk$ is a complex representation $V$
that is equipped with an inner product such that
\[
\pi(A)\quad\mbox{is skew symmetric for all }A\in\gk.
\]
Since $e^{\pi(A)}$ is unitary if and only if $\pi(A)$ is skew symmetric, our
definitions imply that a representation $(V,\Pi)$ of a Lie group $G$ is
unitary if and only if the associated representation $(V,\pi)$ of the real Lie
algebra $\gk$ of $G$ is unitary.

\bt[Compact Lie groups]
Let\label{T:Liecomp} $K$ be a compact Lie group and let $V$ be a
representation of $K$. Then it is possible to equip $V$ with an inner product
so that $V$ becomes a unitary representation of $K$.
\et
\bpro[(sketch)]
Choose an arbitrary inner product $\li\,\cdot\,|\,\cdot\,\re$ on $V$ and
define
\[
\li v|w\re_K:=\int\li\Pi(A)v|\Pi(A)w\re\di A,
\]
where $\di A$ denotes the \emph{Haar measure} on $K$, which is finite by the
assumption that $K$ is compact. It is easy to check that
$\li\,\cdot\,|\,\cdot\,\re_K$ is an inner product. In particular, since
$\Pi(A)$ is invertible for each $A\in K$, we have $\Pi(A)v\neq 0$ and hence
$\li\Pi(A)v|\Pi(A)v\re>0$ for all $v\in V$ and $A\in K$. Now
by the fact that the Haar measure is invariant under the action of the group
\[\ba{l}
\dis\li\Pi(B)v|\Pi(B)w\re_K
=\int\li\Pi(A)\Pi(B)v|\Pi(A)\Pi(B)w\re\di A\\[5pt]
\dis\quad=\int\li\Pi(AB)v|\Pi(AB)w\re\di A=\int\li\Pi(C)v|\Pi(C)w\re\di C
=\li v|w\re_K,
\ea\]
which proves that $V$, equipped with the inner product
$\li\,\cdot\,|\,\cdot\,\re_K$, is a unitary representation of $K$.
\epro

The following lemma is a sort of converse to Theorem~\ref{T:Liecomp} since it
says that noncompact Lie groups do not have faithful unitary representations,
at least when we restrict ourselves to finite-dimensional representations, as
we do here.

\bl[Noncompact Lie groups]
Let\label{L:noncomp} $K$ be a noncompact Lie group and let $V$ be a
faithful (finite dimensional) representation of $K$. Then it is not possible
to equip $V$ with an inner product so that $V$ becomes a unitary
representation of $K$.
\el
\bpro
Equip $V$ with an inner product and let ${\rm U}(V)$ denote the group of all
unitary maps $A:V\to V$. If $(V,\Pi)$ is a faithful representation of $K$,
then the image $\Pi(K)$ of $K$ under $\Pi$ is a closed subset of ${\rm GL}(V)$
and $\Pi:K\to\Pi(K)$ is a homeomorphism. If $(V,\Pi)$ is a unitary
representation, then $\Pi(K)\sub{\rm U}(V)$ and hence by the compactness of
the latter, $\Pi(K)$ is compact. Since $\Pi:K\to\Pi(K)$ is a homeomorphism, it
follows that $K$ is compact.
\epro

A \emph{$\ast$-algebra} is a complex algebra on which there is
defined an \emph{adjoint} operation $A\mapsto A^\ast$ such that
\begin{enumerate}
\item $A\mapsto A^\ast$ is conjugate linear, 
\item $(A^\ast)^\ast=A$,
\item $(AB)^\ast=B^\ast A^\ast$.
\end{enumerate}
If $V$ is a complex finite dimensional linear space equipped with an inner
product, then $\Li(V)$, equipped with the adjoint operation (\ref{adj}), is a
$\ast$-algebra.

A \emph{$\ast$-algebra homomorphism} is an algebra homomorphism
that satisfies
\[
\phi(A^\ast)=\phi(A)^\ast.
\]
A \emph{sub-$\ast$-algebra} of a $\ast$-algebra is a subalgebra that is closed
under the adjoint operation. By definition, a \emph{$\ast$-representation} of
a $\ast$-algebra $\ak$ is a representation $(V,\pi)$ such that $V$ is equipped
with an inner product and $\pi$ is a $\ast$-algebra homomorphism.

In general, a $\ast$-algebra may fail to have a faithful
$\ast$-representation. For finite dimensional $\ast$-algebras, a necessary and
sufficient condition for the existence of a faithful representation is that
\[
A^\ast A=0\quad\mbox{implies}\quad A=0,
\]
but it is rather difficult to prove this; see \cite{Swa17} and references
therein. In infinite dimensions, one needs the theory of C$\ast$-algebras,
which are $\ast$-algebras equipped with a norm that in faithful
representations corresponds to the operator norm $\|A\|=\sup_{\|v\|\leq
  1}\|Av\|$.

Recall the definition of an \emph{adjoint} operation on a complex Lie algebra
$\gk$ from Section~\ref{S:Lieint}. Recall also that we called a Lie algebra
homomorphism \emph{unitary} if $\phi(A^\ast)=\phi(A)^\ast$, and that a
\emph{unitary} representation is a representation $(V,\pi)$ such that $V$ is
equipped with an inner product and $\pi$ is a unitary Lie algebra
homomorphism.

\bl[Universal enveloping $\ast$-algebra]
Let\label{L:Uad} $\gk$ be a Lie-$\ast$-algebra. Then there exists a unique
adjoint operation on its universal enveloping algebra $U(\gk)$ that coincides
with the adjoint operation on $\gk$.
\el
\bpro
Recall from Sections~\ref{S:rep} that every complex linear
space $V$ has a \emph{conjugate} space which is a linear space $\ov V$
together with a conjugate linear bijection $V\ni v\mapsto\ov v\in\ov V$.
If $\ak$ is a complex algebra, then we can equip $\ov\ak$ with the
structure of an algebra by putting
\[
\ov{A}\;\ov{B}:=\ov{BA}.
\]
It is not hard to see that a map $A\mapsto A^\ast$ defined on some algebra
$\ak$ is an adjoint operation if and only if the map $A\mapsto\ov{A^\ast}$
from $\ak$ into $\ov{\ak}$ is an algebra homomorphism. By the definition of an
adjoint operation on a Lie algebra, $[A^\ast,B^\ast]=-[A,B]^\ast$ for all
$A,B\in\gk$. It follows that the map
\[
\gk\ni X\mapsto\ov{X^\ast}\in\ov{U(\gk)}
\]
is a Lie algebra homomorphism, which by the defining property of the universal
enveloping algebra (Theorem~\ref{T:unenv}~(i)) extends to a unique algebra
homomorphism from $U(\gk)$ to $\ov{U(\gk)}$.
\epro

\subsection{Dual, quotient, sum, and product spaces}

\subsubsection*{Dual spaces}

The \emph{dual} $V'$ of a finite dimensional linear space $V$ over $\K=\R$ or
$=\C$ is the space of all linear forms $l:V\to\K$. Each element $v\in V$
naturally defines a linear form $L_v$ on $V'$ by $L_v(l):=l(v)$ and each
linear form on $V$ arises in this way, so we can identify $V''\cong V$.
If $\{e(1),\ldots,e(n)\}$ is a basis for $V$, then setting
$f(i)(e(j)):=1_{\{i=j\}}$ defines a basis $\{f(1),\ldots,f(n)\}$ for $V'$
called the \emph{dual basis}. If $V$ is equipped with an inner product, then
setting
\[
\li v|w:=\li v|w\re
\]
defines a linear form on $V$ and $V':=\{\li v|:v\in V\}$. Through this
identification, we also equip $V'$ with an inner product. Then if
$\{e(1),\ldots,e(n)\}$ is an orthonormal basis for $V$, the dual basis is an 
orthonormal basis for $V'$. Each linear map
$A:V\to W$ gives naturally rise to a \emph{dual map} $A':W'\to V'$ defined by
\[
A'(l):=l\circ A,
\]
and indeed every linear map from $W'$ to $V'$ arises in this way, i.e.,
$\Li(W',V')=\{A':A\in\Li(V,W)\}$. If $V,W$ are equipped with inner products
and $A\in\Li(V,W)$, then
\[
A'(\li\phi|)=\li A^\ast\phi|,
\]
where $A^\ast$ denotes the adjoint of $A$, i.e., this is the linear map
$A^\ast\in\Li(W,V)$ defined by
\[
\li\phi|A\psi\re=\li A^\ast\phi|\psi\re\qquad(\phi\in W,\ \psi\in V).
\]
If $(V,\Pi)$ is a representation of a Lie group $G$, then we can define group
homomorphism $\Pi':G\to\Li(V')$ by
\[
\Pi'(A)l:=\Pi(A^{-1})'l=l\circ\Pi(A^{-1}).
\]
In this way, the dual space $V'$ naturally obtains the structure of a
representation of $G$. Note that
\[
\Pi'(AB)l=l\circ\Pi((AB)^{-1})=l\circ\Pi(A^{-1})\Pi(B^{-1})
=\Pi'(A)(\Pi'(B)l),
\]
proving that $\Pi'$ is indeed a group homomorphism.
Similarly, if $(V,\pi)$ is a representation of a Lie algebra $\gk$, then we can
equip the dual space $V'$ with the structure of a representation of $\gk$ by
putting
\[
\pi'(A)l:=-\pi(A)'(l)=-l\circ\pi(A),
\]
where in this case the minus sign guarantees that
\[\ba{l}
\dis\pi'([A,B])l=-l\circ\pi([A,B])
=-l\circ\big(\pi(A)\pi(B)-\pi(B)\pi(A)\big)\\[5pt]
\dis\quad=-\big(\pi'(B)(\pi'(A)l)-\pi'(A)(\pi'(B)l)
=\pi'(A)(\pi'(B)l)-\pi'(B)(\pi'(A)l).
\ea\]
This is called the \emph{dual representation} or \emph{contragredient
  representation} of $G$ or $\gk$, respectively, associated with $V$, see
\cite[Def.~4.21]{Hal03}. If two representations of $G$ and $\gk$ are
associated as in (\ref{asrep}), then their dual representations are also
associated.

\subsubsection*{Quotient spaces}

By definition, a \emph{normal subgroup} of a group $\Gi$ is a subgroup $\Hi$
such that
\be\label{normal}
A\Hi:=\{AB:B\in\Hi\}=\{BA:B\in\Hi\}=:\Hi A\qquad\forall A\in\Gi,
\ee
or equivalently, if $B\in\Hi$ implies $ABA^{-1}\in\Hi$ for all $A\in\Gi$.
Sets of the form $A\Hi$ and $\Hi A$ are called \emph{left} and \emph{right
cosets}, respectively. If $\Hi$ is a normal subgroup, then left cosets are
right cosets and vice versa, and we can equip the set
\[
\Gi/\Hi:=\big\{A\Hi:A\in\Gi\}=\big\{\Hi A:A\in\Gi\}
\]
of all cosets with a group structure such that
\[
(A\Hi)(B\Hi)=(AB)\Hi.
\]
We call $\Gi/\Hi$ the \emph{quotient group} of $\Gi$ and $\Hi$. Note that as a
set this is obtained from $\Gi$ by dividing out the equivalence relation
\[
A\sim B\quad\desd\quad A=BC\quad\mbox{for some}\quad C\in\Hi.
\]

If $V$ is a linear space and $W\sub V$ is a linear subspace, then we can
define an equivalence relation on $V$ by setting
\[
v_1\sim v_2\quad\desd\quad v_1=v_2+w\quad\mbox{for some}\quad w\in W.
\]
The equivalence classes with respect to this equivalence relation are the sets
of the form
\[
v+W:=\{v+w:w\in W\}
\]
and we can equip the space
\[
V/W:=\{v+W:v\in V\}
\]
with the structure of a linear space by setting
\[
a_1(v_1+W)+a_2(v_2+W):=\big(a_1v_1+a_2v_2\big)+W.
\]

An \emph{invariant subspace} of a representation $V$ of a Lie group $G$, Lie
algebra $\gk$, or algebra $\ak$ is a linear space $W\sub V$ such that $Aw\in
W$ for all $w\in W$ and $A$ from $G$, $\gk$, or $\ak$, respectively. If $W$ is
an invariant subspace, then we can equip the quotient space $V/W$ with the
structure of a representation by setting
\[
A(v+W):=(Av)+W.
\]
Note that this is a good definition since $v_1=v_2+w$ for some $w\in W$
implies $Av_1=Av_2+Aw$ where $Aw\in W$ by the assumption that $W$ is
invariant.

A \emph{left ideal} (resp.\ \emph{right ideal}) of an algebra $\ak$ is a
linear subspace $\ik\sub\ak$ such that $AB\in\ik$ (resp.\ $BA\in\ik$) for all
$A\in\ak$ and $B\in\ik$. An \emph{ideal} is a linear subspace that is both a
left and right ideal. If $\ik$ is an ideal of $\ak$, then we can equip the
quotient space $\ak/\ik$ with the structure of an algebra by putting
\[
(A+\ik)(B+\ik):=(AB)+\ik.
\]
To see that this is a good definition, write $A_1\sim A_2$ if $A_1=A_2+B$ for
some $B\in\ik$. Then $A_1\sim A_2$ and $B_1\sim B_2$ imply that
$A_1=A_2+C$ and $B_1=B_2+D$ for some $C,D\in\ik$ and hence
\[
A_1B_1=(A_2+C)(B_2+D)=A_2B_2+\big(CB_2+A_2D+CD)
\]
with $CB_2+A_2D+CD\in\ik$, so $A_1B_1\sim A_2B_2$. If $\ak$ is a
$\ast$-algebra, then a \emph{$\ast$-ideal} of $\ak$ is an ideal $\ik$ such
that $A\in\ik$ implies $A^\ast\in\ik$. If $\ik$ is a $\ast$-ideal, then we can
equip the quotient algebra $\ak/\ik$ with an adjoint operation by putting
\[
(A+\ik)^\ast:=A^\ast+\ik.
\]

A linear subspace $\hk$ of a Lie algebra $\gk$ is said to be an \emph{ideal}
if $[A,B]\in\hk$ for all $A\in\gk$ and $B\in\hk$. Note that this automatically
implies that also $[B,A]=-[A,B]\in\hk$. If $\hk$ is an ideal of a
Lie algebra, then we can equip the quotient space $\gk/\hk$ with the structure
of a Lie algebra by putting
\[
[A+\hk,B+\hk]:=[A,B]+\hk.
\]
The proof that this is a good definition is the same as for algebras.

\subsubsection*{The direct sum}

The direct sum $V_1\oplus\cdots\oplus V_n$ of linear spaces $V_1,\ldots,V_n$
has already been defined in Section~\ref{S:sumprod}. There is a natural
isomorphism between $V_1\oplus\cdots\oplus V_n$ and the Carthesian product
\[
V_1\times\cdots\times V_n=\big\{\big(\phi(1),\ldots,\phi(n)\big):\phi(i)\in
V_i\ \forall i\big\},
\]
which we equip with a linear structure by defining
\[
a\big(\phi(1),\ldots,\phi(n)\big)
+b\big(\psi(1),\ldots,\psi(n)\big)
:=\big(a\phi(1)+b\phi(1),\ldots,a\phi(n)+b\phi(n)\big).
\]
If $V_1,\ldots,V_n$ are equipped with inner products, then we require that the
inner product on $V_1\oplus\cdots\oplus V_n$ is given by
\be\label{suminp}
\li\phi(1)+\cdots+\phi(n)|\psi(1)+\cdots+\psi(n)\re
:=\sum_{k=1}^n\li\phi(k)|\psi(k)\re,
\ee
which has the effect that $V_1,\ldots,V_n$ are (mutually) orthogonal.
One has the natural isomorphism
\[
(V_1\oplus V_2)/V_2\cong V_1.
\]
In general, given a subspace $V_1$ of some larger linear space $W$, there are
many possible ways to choose another subspace $V_2$ such that $W=V_1\oplus
V_2$ and hence $W\cong(W/V_1)\oplus V_1$.

If $V$ is a linear subspace of some larger linear space $W$, and $W$ is
equipped with an inner product, then we define the \emph{orthogonal
  complement} of $V$ as
\[
V^\perp:=\{w\in W:\li v|w\re=0\ \forall v\in V\}.
\]
Then one has the natural isomorphisms
\[
W/V\cong V^\perp\quand W\cong V\oplus V^\perp,
\]
where the inner product $V\oplus V^\perp$ is given in terms of the inner
products on $V$ and $V^\perp$ as in (\ref{suminp}). Thus, given a linear
subspace $V_1$ of a linear space $W$ that is equipped with an inner product,
there is a canonical way to choose another subspace $V_2$ such that
$W=V_1\oplus V_2$.

If $V_1,\ldots,V_n$ are representations of the same Lie group, Lie algebra, or
algebra, then we equip $V_1\oplus\cdots\oplus V_n$ with the structure of a
representation by putting
\[
A\big(\phi(1)+\cdots+\phi(n)\big):=A\phi(1)+\cdots+A\phi(n).
\]
If $V,W$ are representations, then $W$ is an invariant subspace of $V\oplus W$
and one has the natural isomorphism of representations $(V\oplus W)/W\cong V$.

If $\ak_1,\ldots,\ak_n$ are algebras, then we equip their direct sum
$\ak_1\oplus\cdots\oplus\ak_n$ with the structure of an algebra by putting
\be\label{algsum}
\big(A(1)+\cdots+A(n)\big)\big(B(1)+\cdots+B(n)\big)
:=A(1)B(1)+\cdots+A(n)B(n).
\ee
If $\ak,\bk$ are algebras, then $\bk$ is an ideal of $\ak\oplus\bk$ and one
has the natural isomorphism $(\ak\oplus\bk)/\bk\cong\ak$. Note that $\bk$ is
not a subalgebra of $\ak\oplus\bk$ since $I\not\in\bk$ (unless $\ak=\{0\}$).
For $\ast$-algebras, we also put
\[
\big(A(1)+\cdots+A(n)\big)^\ast:=\big(A(1)^\ast+\cdots+A(n)^\ast\big).
\]
The direct sum of Lie algebras has already been defined in
Section~\ref{S:sumprod}. It is easy to see that this is consistent with the
definition of the direct sum of algebras.

\subsubsection*{The tensor product}

The \emph{tensor product} of two (or more) linear spaces has already been
defined in Section~\ref{S:sumprod}. A proof similar to the proof of
Lemma~\ref{L:covuni} shows that the tensor product is unique up to natural
isomorphisms, i.e., if $V\ti\otimes W$ and
$(\phi,\psi)\mapsto\phi\ti\otimes\psi$ are another linear space and bilinear
map which satisfy the defining property of the tensor product, then there
exists a unique linear bijection $\Psi:V\otimes W\to V\ti\otimes W$ such that
$\Psi(V\otimes W)=V\ti\otimes W$.

If $V,W$ are representations of the same Lie group, then we equip $V\otimes W$
with the structure of a representation by putting
\be\label{Gtens}
A(\phi\otimes\psi):=A\phi\otimes A\psi.
\ee
If $V,W$ are representations of the same Lie algebra or algebra, then we equip
$V\otimes W$ with the structure of a representation by putting
\be\label{altens}
A(\phi\otimes\psi):=A\phi\otimes\psi+\phi\otimes A\psi.
\ee
The reason why we define things in this way is that in view of (\ref{asrep}),
if $\gk$ is the Lie algebra of $G$, then the representation of $\gk$ defined
in (\ref{altens}) is the representation of $\gk$ that is associated with the 
representation of $G$ defined in (\ref{Gtens}). Note that (\ref{altens}) 
is bilinear in $\phi$ and $\psi$ and hence by the defining property of the
tensor product uniquely defines a linear operator on $V\otimes W$.

If $\ak,\bk$ are algebras, then we equip their tensor product $\ak\otimes\bk$
with the structure of an algebra by putting
\[
\big(A(1)\otimes B(1)\big)\big(A(2)\otimes B(2)\big)
:=\big(A(1)A(2)\otimes B(1)B(2)\big).
\]
Using the defining property of the tensor product, one can show that this
unambiguously defines a linear map
\[
(\ak\otimes\bk)^2\ni(A,B)\mapsto AB\in\ak\otimes\bk.
\]
We can identify $\ak$ and $\bk$ with the subalgebras of $\ak\otimes\bk$
given by
\[
\ak\cong\{A\otimes I:A\in\ak\}
\quand
\bk\cong\{I\otimes B:B\in\bk\}.
\]
Note that if we identify $\ak$ and $\bk$ with subalgebras of $\ak\otimes\bk$,
then every element of $\ak$ commutes with every element of $\bk$. If $\ak,\bk$
are $\ast$-algebras, then we equip the algebra $\ak\otimes\bk$ with an adjoint
operation by setting
\[
(A\otimes B)^\ast:=(A^\ast\otimes B^\ast).
\]

If $\gk$ and $\hk$ are Lie algebras, then the universal enveloping algebra of
their direct sum is naturally isomorphic to the tensor product of their
universal enveloping algebras:
\be\label{Uplustim}
U(\gk\oplus\hk)\cong U(\gk)\otimes U(\hk).
\ee
Indeed, if $\{X_1,\ldots,X_n\}$ is a basis for $\gk$ and $\{Y_1,\ldots,Y_m\}$
is a basis for $\hk$, then we can define a bilinear map $(A,B)\mapsto A\otimes
B$ from $U(\gk)\times U(\hk)$ into $U(\gk\oplus\hk)$ by
\[
\ba{l}\dis\big(X_1^{k_1}\cdots X_n^{k_n},Y_1^{l_1}\cdots Y_m^{l_m}\big)\\[5pt]
\dis\quad\mapsto X_1^{k_1}\cdots X_n^{k_n}\otimes Y_1^{l_1}\cdots Y_m^{l_m}
:=X_1^{k_1}\cdots X_n^{k_n}Y_1^{l_1}\cdots Y_m^{l_m}.
\ea\]
where we view $\gk$ and $\hk$ as sub-Lie-algebras of $\gk\oplus\hk$ such that
$[X,Y]=0$ for each $X\in\gk$ and $Y\in\hk$. In view of Theorem~\ref{T:unenv},
the space $U(\gk\oplus\hk)$ together with this bilinear map is a realization
of the tensor product $U(\gk)\otimes U(\hk)$.

On a philosophical note, recall that elements of a Lie algebra are related to
elements of a matrix Lie group via an exponential map. We can view
(\ref{Uplustim}) as a reflection of the property of the exponential map that
converts sums into products.

If $V$ and $W$ are representations of algebras $\ak$ and $\bk$, respectively,
then we can make $V\otimes W$ into a representation of $\ak\otimes\bk$ by
setting
\be\label{AotimB}
(A\otimes B)(\phi\otimes\psi):=(A\phi)\otimes(B\psi).
\ee
Again, by bilinearity and the defining property of the tensor product, this is
a good definition. Note that this is consistent with (\ref{Uplustim}) and our
definition in (\ref{gplush}) where we showed that if $V$ and $W$ are
representations of Lie algebras $\gk$ and $\hk$, then $V\otimes W$ is
naturally a representation of $\gk\oplus\hk$. On the other hand, one should
observe that in the special case that $\ak=\bk$, our present construction
differs from our earlier construction in (\ref{altens}).

\subsection{Irreducible representations}

Let $\gk$ be a Lie algebra on which an adjoint operation is defined, and let
$\hk:=\{\ab\in\gk:\ab^\ast=-\ab\}$ denote the real sub-Lie-algebra\footnote{To
  see that this is a sub-Lie-algebra, note that $\ab,\bb\in\hk$ imply
  $[\ab,\bb]^\ast=-[\ab^\ast,\bb^\ast]$ and hence $[\ab,\bb]\in\hk$.}
consisting of all skew-symmetric elements of $\gk$. It is not hard to see that
$\gk$ is the complexification of $\hk$, i.e., each $\ab\in\gk$ can uniquely be
written as $\ab=\ab_1+i\ab_2$ with $\ab_1,\ab_2\in\hk$.\footnote{Equivalently,
  we may show that each $\ab\in\gk$ can uniquely be written as $\ab={\rm
    Re}(\ab)+i{\rm Im}(\ab)$ with ${\rm Re}(\ab),{\rm Im}(\ab)$
  self-adjoint. This follows easily by putting ${\rm
    Re}(\ab):=\ha(\ab+\ab^\ast)$ and ${\rm Im}(\ab):=\ha i(\ab^\ast-\ab)$.}
Let $\{\xb_1,\ldots,\xb_n\}$ be a basis for $\gk$. The Lie bracket on $\gk$ is
uniquely characterized by the commutation relations
\be\label{comrel2}
[\xb_i,\xb_j]=\sum_{j=1}^nc_{ijk}\xb_k,
\ee
where $c_{ijk}$ are the structure constants (see (\ref{comrel2})). Likewise,
the adjoint operation on $\gk$ is uniquely characterized by its action on
basis elements
\be\label{adrel}
\xb_i^\ast=\sum_jd_{ij}\xb_j,
\ee
where $d_{ij}$ is another set of constants.

By Theorem~\ref{T:Lie3}, the real Lie algebra $\hk$ is the Lie algebra of some
Lie group $G$. By going to the universal cover, we can take $G$ to be simply
connected, in which case it is uniquely determined by $\hk$.
Conversely, if $G$ is a simply connected Lie group, $\hk$ is its real Lie
algebra, and $\gk:=\hk_\C$ is the complexification of $\hk$, then we can equip
$\gk$ with an adjoint operation such that the set of skew symmetric elements
is exactly $\hk$, by putting $(\ab_1+i\ab_2)^\ast:=-\ab_1+i\ab_2$ for each
$\ab_1,\ab_2\in\hk$.

If $V$ is a linear space and $X_1,\ldots,X_n\in\Li(V)$ satisfy (\ref{comrel2}),
then there exists a unique Lie algebra homomorphism $\pi:\gk\to\Li(V)$ such
that $\pi(\xb_i)=X_i$ $(i=1,\ldots,n)$. If $V$ is equipped with an inner
product and the operators $X_1,\ldots,X_n$ moreover satisfy (\ref{adrel}),
then $\pi$ is a unitary representation. By Theorem~\ref{T:unenv}~(i) and
Lemma~\ref{L:Uad}, $\pi$ can in a unique way be extended to a $\ast$-algebra
homomorphism $\ov\pi:U(\gk)\to\Li(V)$. Moreover, if $G$ is the simply connected
Lie group associated with $\hk$, then by Theorem~\ref{T:simpcon}, there
exists a unique Lie group homomorphism $\Pi:G\to\Li(V)$ such that (\ref{Pipi})
holds, so $(V,\Pi)$ is a representation of $G$. Since every element of $\hk$
is skew symmetric, $(V,\pi)$ and hence also $(V,\Pi)$ are unitary
representations of $\hk$ and $G$, respectively.

Let $W\sub V$ be a linear subspace. It is not hard to see that
\[\ba{l}
\dis W\mbox{ is an invariant subspace of }(V,\Pi)\\
\dis\quad\desd\quad
W\mbox{ is an invariant subspace of }(V,\pi)\\
\dis\quad\desd\quad
W\mbox{ is an invariant subspace of }(V,\ov\pi).
\ea\]
We say that $V$ is \emph{irreducible} if its only invariant subspaces are
$\{0\}$ and $V$.

Let $V,W$ be two representations of the same Lie group $G$, Lie algebra $\gk$,
or algebra $\ak$. Generalizing our earlier definition for ie algebras, a
\emph{homomorphism} of representations (of any kind) is a linear map
$\phi:V\to W$ such that
\be\label{twine}
\phi(\ab v)=\ab\phi(v)
\ee
for all $\ab\in G$, $\ab\in\gk$, or $\ab\in\ak$, respectively. Homomorphisms of
representations are called \emph{intertwiners} of representations. If
$\phi$ is a bijection, then its inverse is also an intertwining map. In this
case we call $\phi$ an \emph{isomorphism} and say that the representations are
\emph{equivalent} (or \emph{isomorphic}). If $G$ is a simply connected Lie
group, $\gk$ its associated complexified Lie algebra, and $U(\gk)$ its universal
enveloping algebra, then it is not hard to see that
\[
\mbox{(\ref{twine}) holds }\forall\ab\in G
\ \desd\ \mbox{(\ref{twine}) holds }\forall\ab\in\gk
\ \desd\ \mbox{(\ref{twine}) holds }\forall\ab\in U(\gk).
\]

The following result can be found in, e.g., \cite[Thm~4.29]{Hal03}. In the
special case of complex Lie algebras, we have already stated this in
Proposition~\ref{P:Schur}.

\bp[Schur's lemma]\quad
\begin{itemize}
\item[{\rm\textbf{(a)}}] Let $V$ and $W$ be irreducible representations of a
  Lie group, Lie algebra, or algebra, and let $\phi:V\to W$ be an
  intertwiner. Then either $\phi=0$ or $\phi$ is an isomorphism.
\item[{\rm\textbf{(b)}}] Let $V$ be an irreducible complex representation of a
  Lie group, Lie algebra, or algebra, and let $\phi:V\to V$ be an
  intertwiner. Then $\phi=\la I$ for some $\la\in\C$.
\end{itemize}
\ep

By definition, the \emph{center} of an algebra is the subalgebra
$\Ci(\ak):=\{C\in\ak:[A,C]=0\ \forall A\in\ak\}$. The center is \emph{trivial}
if $\Ci(\ak)=\{\la I:\la\in\K\}$. The following is adapted from
\cite[Cor.~4.30]{Hal03}.

\bcor[Center]
Let\label{C:cent} $(V,\pi)$ be an irreducible complex representation of an
algebra $\ak$ and let $C\in\Ci(\ak)$. Then $\pi(C)=\la I$ for some $\la\in\C$.
\ecor
\bpro
Define $\phi:V\to V$ by $\phi v:=\pi(C)v$. Then
$\phi(Av)=\pi(C)\pi(A)v=\pi(CA)v=\pi(AC)v=\pi(A)\pi(C)v=A(\phi v)$ for all
$A\in\ak$, so $\phi:V\to V$ is an intertwiner. By part~(b) of Schur's lemma,
$\phi=\la I$ for some $\la\in\C$.
\epro

\subsection{Semisimple Lie algebras}

A Lie algebra $\gk$ is called \emph{irreducible} (see
\cite[Def.~3.11]{Hal03}) if its only ideals are $\{0\}$ and $\gk$, and
\emph{simple} if it is irreducible and has dimension ${\rm dim}(\gk)\geq 2$.
A Lie algebra is called \emph{semisimple} if it can be written as the direct
sum of simple Lie algebras. Recall the definition of the center of a Lie
algebra in~(\ref{center}).

\bl[Trivial center]
The\label{L:trivcent} center of a semisimple Lie algebra is trivial.
\el
\bpro
If $\gk$ is simple and $A$ is an element of its center, then the linear space
spanned by $A$ is an ideal. Since ${\rm dim}(\gk)\geq 2$ and its only ideals
are $\{0\}$ and $\gk$, this implies that $A=0$. If
$\gk=\gk_1\oplus\cdots\oplus\gk_n$ is the direct sum of simple Lie algebras,
then we can write any element $A$ of the center of $\gk$ as
$A=A(1)+\cdots+A(n)$ with $A(k)\in\gk$. By the definition of the Lie bracket
on $\gk$ (see (\ref{Liesum})), $A(k)$ lies in the center of $\gk$ for each
$k$, and hence $A=0$ by what we have already proved.
\epro

The following proposition is similar to \cite[Prop.~7.4]{Hal03}.

\bp[Inner product on Lie algebra]
Let\label{P:gkin} $\gk$ be a Lie algebra on which an adjoint operation is defined, let
$\hk:=\{\ab\in\gk:\ab^\ast=-\ab\}$, and let $G$ be the simply connected Lie
group with Lie algebra $\hk$. Assume that $G$ is compact. Then the Lie algebra
$\gk$, equipped with the map
\[
\gk\ni\xb\mapsto{\rm ad}_\xb\in\Li(\gk),
\]
is a faithful representation of itself. It is possible to equip $\gk$ with an
inner product such that this is a unitary representation, i.e., ${\rm
  ad}_{\xb^\ast}=({\rm ad}_\xb)^\ast$ $(\xb\in\gk)$.
\ep
\bpro
By \cite[Prop.~7.7]{Hal03}, the center of $\gk$ is trivial. By
Lemma~\ref{L:Lierep} and the remarks below it, it follows that $\gk$, equipped
with the map $\gk\ni{\rm ad}_X\in\Li(\gk)$, is a faithful representation of
itself.  This representation naturally gives rise to a representation of $G$.
By assumption, $G$ is compact, so by Theorem~\ref{T:Liecomp}, we can equip
$\gk$ with an inner product so that this representation is unitary. It follows
that the representation of $\hk$ on $\gk$ is also unitary and hence the
representation of $\gk$ on itself is a unitary representation.
\epro

The following theorem follows from \cite[Thm~7.8]{Hal03}.

\bt[Semisimple algebras]
Let\label{T:comp} $G$ be a compact simply connected Lie group and let $\gk$ be
the complexification of its Lie algebra. Then $\gk$ is semisimple.
\et
\bpro[(main idea)]
If $\gk$ is not simple, then it has some ideal $\ik$ that is neither $\{0\}$
nor $\gk$. Let $\ik^\perp$ denote the orthogonal complement of $\ik$ with
respect to the inner product on $\gk$ defined in Proposition~\ref{P:gkin}. It is
shown in \cite[Prop.~7.5]{Hal03} that $\ik^\perp$ is an ideal of $\gk$ and one
has $\gk\cong\ik\oplus\ik^\perp$, where $\oplus$ denotes the direct sum of Lie
algebras. Continuing this process, one arrives at a decomposition of $\gk$ as
a direct sum of simple Lie algebras.
\epro

In fact, the converse statement to Theorem~\ref{T:comp} also holds: if $\gk$
is a semisimple complex Lie algebra, then it is the complexification of the
Lie algebra of a compact simply connected Lie group. This is stated (with
references for a proof) in \cite[Sect.~10.3]{Hal03}.

Let $G$ be a compact simply connected Lie group, let $\hk$ be its real Lie
algebra, let $\gk:=\hk_\C$ be the complexification of $\hk$, and let $U(\gk)$
denote the universal enveloping algebra of $\gk$. The \emph{Casimir element}
is the element $C\in U(\gk)$ defined as
\[
\cb:=-\sum_j\xb_j^2,
\]
where $\{\xb_1,\ldots,\xb_n\}$ is a basis for $\hk$ that is orthonormal with
respect to the inner product from Proposition~\ref{P:gkin}.\footnote{The inner
  product from Proposition~\ref{P:gkin} is not completely unique; at best it
  is only determined up to a multiplicative constant. So the
  Casimir operator depends on the choice of the inner product, but once this
  is fixed, it does not depend on the choice of the orthonormal basis.} We cite
the following result from \cite[Prop.~10.5]{Hal03}.

\bp[Casimir element]
The\label{P:Casimir} definition of the Casimir element does not depend on the
choice of the orthonormal basis $\{\xb_1,\ldots,\xb_n\}$ of $\hk$. Moreover
$\cb$ lies in the center of $U(\gk)$.
\ep

In irreducible representations, the Casimir element has a simple form.

\bl[Representations of Casimir element]
For\label{L:Caspos} each irreducible representation $(V,\pi)$ of $\gk$, there
exists a constant $\la_V\geq 0$ such that $\pi(\cb)=\la_VI$.
\el
\bpro
Proposition~\ref{P:Casimir} and Corollary~\ref{C:cent} imply that for each
irreducible representation $(V,\pi)$ of $U(\gk)$, there exists a constant
$\la\in\C$ such that $\pi(\cb)=\la I$. By Theorem~\ref{T:Liecomp}, we can equip
$V$ with an inner product such that it is a unitary representation of $\hk$.
This means that $\xb_j$ is skew symmetric and hence $i\xb_j$ is hermitian, so
$\cb=\sum_i(i\xb_j)^2$ is a positive operator. In particular, its eigenvalues
are $\geq 0$.
\epro

\subsection{Some basic matrix Lie groups}

For any finite-dimensional linear space $V$ over $V=\R$ or $=\C$, we let ${\rm
  GL}(V)$ denote the \emph{general linear group} of all invertible linear maps
$A:V\to V$. In particular, we write ${\rm GL}(n;\R)={\rm GL}(\R^n)$ and ${\rm
  GL}(n;\C)={\rm GL}(\C^n)$.

The \emph{special linear group} ${\rm SL}(V)$ is defined as
\[
{\rm SL}(V):=\big\{A\in{\rm GL}(V):{\rm det}(A)=1\big\}.
\]
Again, we write ${\rm SL}(n;\R)={\rm SL}(\R^n)$ and ${\rm SL}(n;\C)={\rm
  SL}(\C^n)$. If $V$ is a finite-dimensional linear space over $\C$ and $V$ is
equipped with an inner product $\li\,\cdot|\,\cdot\,\re$, then we call
\[
{\rm U}(V):=\{A\in\Li(V):A\mbox{ is unitary}\}
\]
the \emph{unitary group} and
\[
{\rm SU}(V):=\{A\in{\rm U}(V):{\rm det}(A)=1\}
\]
the \emph{special unitary group}, and write ${\rm U}(n):={\rm U}(\C^n)$
and ${\rm SU}(n):={\rm SU}(\C^n)$.

If $V$ is a finite-dimensional linear space over $\R$ and $V$ is
equipped with an inner product $\li\,\cdot|\,\cdot\,\re$, then an operator
$O\in\Li(V)$ that preserves the inner product as in (\ref{Udef}) is called
\emph{orthogonal}. (This is the equivalent of unitarity in the real setting.)
We call
\[
{\rm O}(V):=\{A\in\Li(V):A\mbox{ is orthogonal}\}
\]
denote the \emph{orthogonal group} and
\[
{\rm SO}(V):=\{A\in{\rm O}(V):{\rm det}(A)=1\}
\]
the \emph{special orthogonal group}, and write ${\rm O}(n):={\rm O}(\R^n)$
and ${\rm SO}(n):={\rm SO}(\R^n)$. There also exists a group ${\rm O}(n;\C)$,
which consists of all complex matrices that preserve the bilinear form
$(v,w):=\sum_iv_iw_i$. Not that this is \emph{not} the inner product on
$\C^n$; as a result ${\rm O}(n;\C)$ is not the same as ${\rm U}(n)$.

Unitary operators satisfy $|{\rm det}(A)|=1$ and orthogonal operators satisfy
${\rm det}(A)=\pm 1$. The group ${\rm O}(3)$ consists of rotations and
reflections (and combinations thereof) while ${\rm SO}(3)$ consists only of
rotations.

By \cite[Prop.~3.23]{Hal03}, for $\K=\R$ or $=\C$, the Lie algebra of ${\rm
  GL}(n,\K)$ is the space $M_n(\K)$ of all $\K$-valued $n\times n$ matrices,
and the Lie algebra of ${\rm SL}(n,\K)$ is given by
\[
\slk(n,\K)=\{A\in M_n(\K):\tr(A)=0\}.
\]
By \cite[Prop.~3.24]{Hal03}, the Lie algebras of ${\rm U}(n)$ and
${\rm O}(n)$ are given by
\[
\uk(n)=\{A\in M_n(\C):A^\ast=-A\}
\quand
\ok(n)=\{A\in M_n(\R):A^\ast=-A\}.
\]
Moreover, again by \cite[Prop.~3.24]{Hal03}, the Lie algebras of ${\rm SU}(n)$
and ${\rm SO}(n)$ are given by
\[
\suk(n)=\{A\in M_n(\C):A^\ast=-A,\ \tr(A)=0\}
\quand
\sok(n)=\ok(n).
\]
By \cite[formula (3.17)]{Hal03}, the complexifications of the real Lie
algebras introduced above are given by
\[\bac
\dis\glk(n,\R)_\C&\cong&\dis\glk(n,\C),\\[5pt]
\dis\uk(n)_\C&\cong&\dis\glk(n,\C),\\[5pt]
\dis\suk(n)_\C&\cong&\dis\slk(n,\C),\\[5pt]
\dis\slk(n,\R)_\C&\cong&\dis\slk(n,\C),\\[5pt]
\dis\sok(n,\R)_\C&\cong&\dis\sok(n,\C).
\ea\]

As mentioned in \cite[Sect.~1.3.1]{Hal03}, the following Lie groups are
compact:
\[
{\rm O}(n),\quad{\rm SO}(n),\quad{\rm U}(n),\quand{\rm SU}(n).
\]
By \cite[Prop~1.11, 1.12, and 1.13]{Hal03} and \cite[Exercise~1.13]{Hal03},
the following Lie groups are connected:
\[
{\rm GL}(n;\C)\quad{\rm SL}(n;\C)\quad{\rm U}(n)\quad{\rm SU}(n),
\quand{\rm SO}(n).
\]
By \cite[Prop.~13.11]{Hal03}, the group ${\rm SU}(n)$ is simply connected.
By \cite[Example~5.15]{Hal03}, ${\rm SU}(2)$ is the universal cover of ${\rm
  SO}(3)$.

Of further interest are the real and complex \emph{symplectic groups}
${\rm SP}(n,\R)$ and ${\rm SP}(n,\C)$, and the \emph{compact symplectic group}
${\rm SP}(n)$; for their definitions we refer to \cite[Sect.~1.2.4]{Hal03}.

\subsection{The Lie group SU(1,1)}

Let us define a Minkowski form $\{\,\cdot\,,\,\cdot\,\}:\C^2\to\C$ by
\[
\{v,w\}:=v_1^\ast w_1-v_2^\ast w_2.
\]
Note that this is almost identical to the usual definition of the inner
product on $\C^2$ (in particular, it is conjugate linear in its first argument and
linear in its second argument), except for the minus sign in front of the
second term. Letting
\[
M:=\left(\ba{cc}1&0\\ 0&-1\ea\right),
\]
we can write
\[
\{v,w\}=\li v|M|w\re,
\]
where $\li\,\cdot\,,\,\cdot\,\re$ is the usual inner product.
The Lie group ${\rm SU}(1,1)$ is the matrix Lie group consisting of all
matrices $Y\in\Li(\C^2)$ with determinant 1 that preserve this Minkowski form,
i.e.,
\[
{\rm det}(Y)=1\quand\{Yv,Yw\}=\{v,w\}\quad(v,w\in\C^2).
\]
The second condition can be rewritten as $\li Yv|M|Yw\re=\li
v|M|w\re$ which holds for all $v,w$ if and only if
\be\label{Munit}
Y^\ast MY=M,
\ee
where $Y^\ast$ denotes the usual adjoint of a matrix. Since
\[
(e^{tA})^\ast Me^{tA}=M+t(A^\ast M+MA)+O(t^2),
\]
it is not hard to see that a matrix of the form $Y=e^{tA}$ satisfies
(\ref{Munit}) if and only if
\[
A^\ast M+MA=0\quad\desd\quad MA^\ast M=-A,
\]
and the Lie algebra $\suk(1,1)$ associated with ${\rm SU}(1,1)$ is given by
\[
\suk(1,1)=\big\{A\in M_2(\C):MA^\ast M=-A,\, \tr(A)=0\big\}.
\]
It is easy to see that
\[
A=\left(\ba{cc}A_{11}&A_{12}\\ A_{21}&A_{22}\ea\right)
\quad\volgt\quad
MA^\ast M=\left(\ba{cc}A_{11}&-(A_{21})^\ast\\-(A_{12})^\ast&A_{22}\ea\right)
\]
and in fact the map $A\mapsto MA^\ast M$ satisfies the axioms of
an adjoint operation. Let $\suk(1,1)_\C$ denote the Lie algebra
\[
\suk(1,1)_\C:=\big\{A\in M_2(\C):\tr(A)=0\big\},
\]
equipped with the adjoint operation $A\mapsto MA^\ast M$. Then $\suk(1,1)$ is
the real sub-Lie algebra of $\suk(1,1)_\C$ consisting of all elements that are
skew symmetric with respect to the adjoint operation $A\mapsto MA^\ast M$.

A basis for $\suk(1,1)_\C$ is formed by the matrices in (\ref{psPauli}),
which satisfy the commutation relations (\ref{su11com}). The adjoint
operation $A\mapsto MA^\ast M$ leads to the adjoint relations (\ref{su11adj}).
Some elementary facts about the Lie algebra $\suk(1,1)_\C$ are already stated
in Section~\ref{S:SU11}. Note that the definition of the ``Casimir operator''
in (\ref{CasimT}) does not follow the general definition for compact Lie
groups in Proposition~\ref{P:Casimir}, but is instead defined in an analogous
way, replacing the inner product by a Minkowski form.

\subsection{The Heisenberg group}

Consider the matrices
\[
X:=\left(\ba{ccc}0&1&0\\ 0&0&0\\ 0&0&0\ea\right),\quad
Y:=\left(\ba{ccc}0&0&0\\ 0&0&1\\ 0&0&0\ea\right),\quad
Z:=\left(\ba{ccc}0&0&1\\ 0&0&0\\ 0&0&0\ea\right).
\]
We observe that
\[\ba{lll}
XX=0,\quad&XY=Z,\quad&XZ=0,\\
YX=0,\quad&YY=0,\quad&YZ=0,\\
ZX=0,\quad&ZY=0,\quad&ZZ=0.
\ea\]
The \emph{Heisenberg group} $H$ \cite[Sect.~1.2.6]{Hal03} is the matrix Lie
group consisting of all $3\times 3$ real matrices of the form
\[
B=I+xX+yY+zZ\qquad(x,y,z\in\R).
\]
To see that this is really a group, we note that if $B$ is
as above, then its inverse $B^{-1}$ is given by
\[
B^{-1}=-xX-yY+(xy-z)Z.
\]
It is easy to see that $\{X,Y,Z\}$ is a basis for the Lie algebra $\hk$ of
$H$. In fact, the expansion formula for $e^{t(xX+yY+zZ)}$ terminates and
\[
e^{t(xX+yY+zZ)}=I+t(xX+yY+zZ)+\ha t^2xyZ\qquad(t\geq 0).
\]
The basis elements $X,Y,Z$ satisfy the commutation relations
\[
[X,Y]=Z,\quad[X,Z]=0,\quad[Y,Z]=0.
\]
Thus, we can abstractly define the \emph{Heisenberg Lie algebra} as the real
Lie algebra $\hk$ with basis elements $\xb,\yb,\zb$ that satisfy the
commutation relations
\be\label{heis}
[\xb,\yb]=\zb,\quad[\xb,\zb]=0,\quad[\yb,\zb]=0.
\ee
Representations of the Heisenberg algebra have already been discussed in
Subsection~\ref{S:Heis}.


\begin{thebibliography}{CGGR13}

\bibitem[AH07]{AH07}
R.~Alkemper and M.~Hutzenthaler.
Graphical representation of some duality relations in stochastic population
models.
\emph{Electron.\ Commun.\ Probab.}~12 (2007), 206--220. 

\bibitem[Bar47]{Bar47}
V.~Bargmann.
Irreducible unitary representations of the Lorentz group. 
\emph{Ann.\ Math.}~48 (1947), 568--640.

\bibitem[Bar61]{Bar61}
V.~Bargmann.
On a Hilbert space of analytic functions and an associated integral transform.
\emph{Commun.\ Pure Appl.\ Math.}~14 (1961), 187--214.

\bibitem[CGGR15]{CGGR15}
G.~Carinci, C.~Giardin\`a, C.~Giberti, and F.~Redig.
Dualities in population genetics: a fresh look with new dualities.
\emph{Stochastic Processes Appl.}~125(3) (2015), 941--969.

\bibitem[Fil92]{Fil92}
J.A.~Fill.
Strong stationary duality for continuous-time Markov chains. I. Theory.
{\em J.\ Theor.\ Probab.}~5(1) (1992), 45--70.

\bibitem[GKRV09]{GKRV09}
C.~Giardin\`a, J.~Kurchan, F.~Redig, and K.~Vafayi.
Duality and hidden symmetries in interacting particle systems.
\emph{J.\ Stat.\ Phys.}~135(1) (2009), 25--55.

\bibitem[GRV10]{GRV10}
C.~Giardin\`a, F.~Redig, and K.~Vafayi.
Correlation Inequalities for Interacting Particle Systems with Duality.
\emph{J.\ Stat.\ Phys.}~141(2) (2010), 242--263.

\bibitem[Hal03]{Hal03}
B.C.~Hall.
\emph{Lie Groups, Lie Algebras, and Representations: An Elementary
  Introduction.}
Graduate Texts in Mathematics, vol. 222, Springer, 2003.

\bibitem[JK14]{JK14}
S.~Jansen and N.~Kurt.
On the notion(s) of duality for Markov processes.
\emph{Prob.\ Surveys}~11 (2014), 59--120.

\bibitem[LS95]{LS95}
A.~Sudbury and P.~Lloyd.
Quantum operators in classical probability theory.~II:
The concept of duality in interacting particle systems.
\emph{Ann.\ Probab.}~23(4) (1995), 1816--1830.

\bibitem[LS97]{LS97}
A.~Sudbury and P.~Lloyd.
Quantum operators in classical probability theory.~IV:
Quasi-duality and thinnings of interacting particle systems.
\emph{Ann.\ Probab.}~25(1) (1997), 96--114.

\bibitem[Nov04]{Nov04}
M.~Novaes.
Some basics of $\suk(1,1)$.
\emph{Revista Brasileira de Ensino de F\'isica}~26(4) (2014), 351--357.

\bibitem[Ros04]{Ros04}
J.~Rosenberg.
A selective history of the Stone-von Neumann theorem.
\emph{Contemp.\ Math.}~365 (2004), 331--353.

\bibitem[SS16]{SS16}
A.~Sturm and J.M.~Swart.
Pathwise duals of monotone and additive Markov processes.
\emph{J.\ Theor.\ Probab.} (2016). 52 pages. doi:10.1007/s10959-016-0721-5.

\bibitem[Sud00]{Sud00}
A.~Sudbury.
Dual families of interacting particle systems on graphs.
\emph{J.\ Theor.\ Probab.}~13(3) (2000), 695--716.

\bibitem[Swa06]{Swa06}
J.M.~Swart.
Duals and thinnings of some relatives of the contact process.
Pages 203-214 in: \emph{Prague Stochastics 2006},
M. Hu\v{s}kov\'a and M. Jan\v{z}ura (eds.), Matfyzpress, Prague, 2006. 
ArXiv:math.PR/0604335.

\bibitem[Swa13]{Swa13}
J.M.~Swart.
\emph{Duality and Intertwining of Markov Chains.}
Lecture notes for the ALEA in Europe School October 21-25 2013, Luminy
(Marseille).
Available at:\\
{\tt http://staff.utia.cas.cz/swart/tea\underline{\ }index.html}.

\bibitem[Swa17]{Swa17}
J.M.~Swart.
\emph{Introduction to Quantum Probability}.
Lecture notes (2017) available at:\\
{\tt http://staff.utia.cas.cz/swart/tea\underline{\ }index.html}.

\bibitem[VK91]{VK91}
N.Ya.~Vilenkin and A.U.~Klimyk. 
\emph{Representation of Lie groups and special functions.
Volume 1: Simplest Lie groups, special functions and integral transforms.}
%Transl.\ from the Russian by V.A.~Groza and A.A.~Groza.
%Mathematics and Its Applications. Soviet Series.~72.
Kluwer, Dordrecht, 1991.

\end{thebibliography}
\end{document}